\RequirePackage{fix-cm}
\documentclass[twocolumn]{svjour3}          
\smartqed  
\usepackage{graphicx}
\usepackage{amsmath, amssymb}
\usepackage{algorithm}
\usepackage[noend]{algpseudocode}
\usepackage{tikz}
\usepackage{multirow}
\usepackage{array}
\usepackage{makecell} 
\usepackage{pifont}
\usetikzlibrary{arrows,decorations.pathreplacing}
\usetikzlibrary{arrows.meta}
\usepackage{cite}
\usepackage[utf8]{inputenc}
\newcommand{\xmark}{\ding{55}}%

\newcommand{\revision}[1]{{\color{black} #1}}

\renewcommand{\vec}[1]{\mathbf{#1}}
\newcommand{\mesh}{\Omega}
\newcommand{\timeint}{(t_0,t_f]}

\newcolumntype{L}[1]{>{\raggedright\arraybackslash}p{#1}}
\newcolumntype{C}[1]{>{\centering\arraybackslash}p{#1}}
\newcolumntype{R}[1]{>{\raggedleft\arraybackslash}p{#1}}

\begin{document}

\title{Parallel-in-Time Simulation of an Electrical Machine using MGRIT
}


\author{Matthias Bolten        \and
        Stephanie Friedhoff \and
        Jens Hahne \and
		Sebastian Sch{\"{o}}ps
}


\institute{Matthias Bolten \and Stephanie Friedhoff \and Jens Hahne \at
              Fakult{\"{a}}t f{\"{u}}r Mathematik und Naturwissenschaften \\
              Bergische Universit{\"{a}}t Wuppertal, 42097 Wuppertal, Germany \\
              \email{bolten@math.uni-wuppertal.de} \\
              \email{friedhoff@math.uni-wuppertal.de}  \\
              \email{jens.hahne@math.uni-wuppertal.de}           
           \and
           Sebastian Sch{\"{o}}ps \at
			Centre for Computational Engineering, Technische Universit{\"{a}}t Darmstadt,
			64293 Darmstadt \\
			\email{schoeps@temf.tu-darmstadt.de} \\
}

\date{Received: date / Accepted: date}

\maketitle

\begin{abstract}
We apply the multigrid-reduction-in-time \linebreak (MGRIT) algorithm to an eddy current simulation of a two-dimensional induction machine supplied by a pulse-width-modulation signal. To resolve the fast-switching excitations, small time steps are needed, such that parallelization in time becomes highly relevant for reducing the simulation time. The MGRIT algorithm is an iterative method that allows calculating multiple time steps simultaneously by using a time-grid hierarchy. It is particularly well suited for introducing time parallelism in the simulation of electrical machines using existing application codes, as MGRIT is a non-intrusive approach that essentially uses the same time integrator as a traditional time-stepping algorithm. However, the key difficulty when using time-stepping routines of existing application codes for the MGRIT algorithm is that the cost of the time integrator on coarse time grids must be less expensive than on the fine grid to allow for speedup over sequential time stepping on the fine grid. To overcome this difficulty, we consider reducing the costs of the coarse-level problems by adding spatial coarsening. We investigate effects of spatial coarsening on MGRIT convergence when applied to two numerical models of an induction machine, one with linear material laws and a full nonlinear model. Parallel results demonstrate significant speedup in the simulation time compared to sequential time stepping, even for moderate numbers of processors.
\end{abstract}

\section{Introduction}
\label{intro}

Induction motors are electrical machines in which the magnetic field in the rotor is obtained by an asynchronous motion with respect to the field in the stator, e.\,g. \cite{Pyrhonen_2008aa}. 
In the design process of such machines, motion and electromagnetic fields are numerically determined by space and time discretization of Maxwell's equations. The spatial discretization, typically by finite elements, may lead to large models, e.\,g., due to resolving the so-called skin effect \cite{jackson_classical_1999}, but this is often mitigated by considering only two-dimensional models. More computational burden comes from the fact that large time intervals must be considered for simulating the start-up phase, i.\,e., until the machine reaches its steady state. Moreover, machines are commonly excited by pulse-width modulated (PWM) excitations. PWM uses fast-switching pulses whose width is controlled such that the time-average corresponds to a specific waveform, typically a sine wave such that very small time steps are needed. As a consequence, the computational complexity of transient simulations is extremely high.

Carrying out the simulations on modern parallel compute systems using existing application codes may still result in runtimes that are beyond what is feasible in practice, since parallelism in existing codes is limited to the spatial aspect of the problem. Due to the stagnation of processor clock speeds in recent years, modern compute systems are characterized by growing numbers of processors. Existing application codes cannot take full advantage of these systems, since spatial parallelism limits the amount of processors that can be exploited. A promising approach for reducing simulation times are parallel-in-time methods that introduce parallelism in the temporal dimension. Especially for simulations over long time intervals with small time steps as considered for electrical machines, parallelization in time can become even more relevant than in space.

To this end, various approaches have been proposed for speeding up the transient simulation of electrical machines, e.\,g., the time-periodic explicit error correction method \cite{Takahashi_2015aa,Takahashi_2019aa} or the extraction of a circuit model \cite{Bermudez_2019aa}.  Also, two of the most well-known parallel-in-time methods, Parareal \cite{JLLions_etal_2001a} and MGRIT \cite{RDFalgout_etal_2014}, have been applied to eddy current problems. In particular, Parareal has been considered for RL-circuit models and an induction machine \cite{Schops_2018aa,Gander_2019aa,Bast_2019aa}, and MGRIT has been applied to a coaxial cable model \cite{friedhoff2019multigridreductionintime,Friedhoff_2020aa}. This paper focuses on MGRIT applied to two numerical models of an induction machine, providing both an extension of ideas of the work on Parareal for problems involving discontinuous or multirate excitations to a multilevel setting as well as broadening the application areas of MGRIT to include more complex eddy current problems.

The idea of MGRIT is to use a time-grid hierarchy, in which each grid level is characterized by a time integrator, with computational costs decreasing from fine to coarse levels. Expensive time-integration routines are applied in parallel on temporal subdomains, and on the coarsest level, time integration is performed over the whole temporal domain. The key difficulty for achieving a significant reduction in simulation time over time-stepping algorithms is the choice of the time integrators such that the cost of the time integrators on coarse time grids is less expensive than on the fine grid. A common approach for defining a time-grid hierarchy is to use different temporal resolutions \cite{RDFalgout_etal_2014,FalgoutEtAl2017b,AHessenthaler_etal_2018}. Other techniques such as considering time-integration routines of different orders of accuracy have also been applied, e.\,g., in power-grid simulations \cite{Lecouvez2016_powersystem,Schroder2018_powersystems}. When using time-stepping routines of existing application codes, the combination of various temporal and spatial resolutions is an attractive approach for defining a time-grid hierarchy, since this allows adding time parallelism non-intrusively and reducing costs on coarse levels at the same time. 

The use of coarser spatial meshes has been explored for various parallel-in-time methods. For Parareal, this approach has been considered in \cite{Ruprecht2014_GAMM,LunetEtAl2018}. Furthermore, the combination of temporal coarsening with coarsening in space is natural in space-time multigrid methods, including the methods in \cite{Hackbusch1984,VandewalleHorton1994,NeumuellerGander2016,FrancoEtAl2018}. The composition of the ``parallel full approximation scheme in space and time'' (PFASST) \cite{EmmettMinion2012} with spatial multigrid, including spatial coarsening, has been considered in \cite{MinionEtAl2015}, scaling results can be found in \cite{SpeckEtAl2014_Parco,RuprechtEtAl2013_SC}. The use of spatial coarsening in MGRIT has been explored for the $p$-Laplacian \cite{MR3716560} and for Burgers equation \cite{Howse2019_burgers}. In this paper, we demonstrate that spatial coarsening can significantly reduce the runtime of numerical simulations of an induction machine, but it may also degrade MGRIT convergence, as observed for the $p$-Laplacian in \cite{MR3716560}. Our tests show that convergence can be improved by using a stronger relaxation scheme.

The remainder of this paper is structured as follows. In Section~\ref{sec:mgrit}, the MGRIT algorithm is introduced for ordinary differential equations, Section~\ref{sec:sim_induction_machine} discusses the induction machine model and its spatial discretization. Numerical results of the parallel-in-time simulation are presented in Sections~\ref{sec:linear_model} and \ref{sec:nonlinear_model}. Finally, conclusions are given in \ref{sec:conclusions}.

\section{Multigrid-Reduction-in-Time}
\label{sec:mgrit}
Consider a system of ordinary differential equations (ODEs) of the form
    \begin{equation}\label{eq:problem_system}
    	\vec{u}'(t) = \vec{f}(t,\vec{u}(t)), \quad \vec{u}(t_0) = \mathbf{g}_0, \quad t \in (t_0,t_f],
    \end{equation}
\revision{with arbitrary time points $t_0$ and $t_f$, where $t_f>t_0$, }arising, for example, after the spatial discretization of a space-time partial differential equation (PDE). Discretize the time interval on a uniform distributed temporal mesh $t_i = i\Delta t, \; i=0,1,\dots,N_t$, with constant time step $\Delta t = (t_f-t_0)/N_t$, and let $\vec{u}_i \approx \vec{u}(t_i)$ for $i = 0,\dots,N_t$ with $\vec{u}_0 = \vec{u}(0)$. Consider a one-step time integration method
	\begin{equation}\label{eq:mgrit_nonlinear_system}
		\vec{u}_i = \boldsymbol{\Phi}_{i} (\vec{u}_{i-1}), \quad i = 1,2,\dots,N_t,
	\end{equation}
with time integrator $\boldsymbol{\Phi}_{i}$, propagating the solution $\vec{u}_{i-1}$ from time point $t_{i-1}$ to time point $t_i$, including forcing from the right-hand-side of the PDE. Abusing notation for the nonlinear case, Equation \eqref{eq:mgrit_nonlinear_system} can be written in matrix form as
	\[
		\revision{\mathcal{A}}\vec{u}\equiv\begin{bmatrix}
			I\\
			-\boldsymbol{\Phi}_{1} & I\\
			& \ddots & \ddots\\
			& & -\boldsymbol{\Phi}_{N_t} & I
		\end{bmatrix}\begin{bmatrix}
			\vec{u}_0\\
			\vec{u}_1\\
			\vdots\\
			\vec{u}_{N_t}
		\end{bmatrix} = \begin{bmatrix}
			\mathbf{g}_0\\
			\mathbf{0}\\
			\vdots\\
			\mathbf{0}
		\end{bmatrix}\equiv \mathbf{g}
	\]
and solved by a sequential forward block solve. This method is optimal, i.\,e., $\mathcal{O}\left(N_t\right)$, and gives the discrete solution after $N_t$ applications of the time integrator. 

\revision{In contrast to this sequential-in-time approach, \linebreak parallel-in-time methods allow for parallelism in the time domain. In the following, we present one of these methods, the multigrid-reduction-in-time (MGRIT) algorithm.}

\subsection{The MGRIT algorithm}
\label{subsec:mgrit-alg}

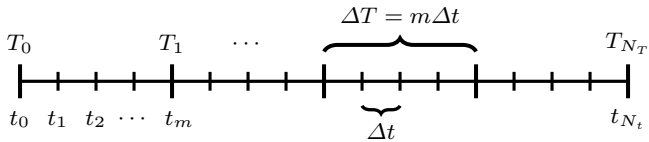
\begin{figure}
\begin{center}
	\begin{tikzpicture}
	\draw[line width=1.2pt] (0,1.3) -- (8,1.3);
	\foreach \i in {0,.5,1,...,8}{
		\draw[line width=1.15pt] (\i,.125+1.3) -- (\i,-.125+1.3);
	}
	\foreach \i in {0,2,...,8}{
		\draw[line width=1.5pt] (\i,.225+1.3) -- (\i,-.225+1.3);
	}
	\foreach \i/\n in {0/0,.5/1,1/2}{
		\draw (\i,-.5+1.3) node {$t_\n$};
	}
	\draw (1.5,-.5+1.3) node {$\cdots$};
	\draw (2.1,-.5+1.3) node {$t_m$};
	\draw (8,-.5+1.3) node {$t_{N_t}$};
	\draw (0,.5+1.3) node {$T_0$};
	\draw (2,.5+1.3) node {$T_1$};
	\draw (3,.5+1.3) node {$\cdots$};
	\draw (8,.5+1.3) node {$T_{N_T}$};
			
	\draw[line width=1.25pt,decorate,decoration={brace,amplitude=3pt},yshift=5pt] (5,-.5+1.3) -- (4.5,-.5+1.3) node[below=3pt,midway] {$\Delta t$};
	\draw[line width=1.5pt,decorate,decoration={brace,amplitude=5pt},yshift=-3pt] (4,.5+1.3) -- (6,.5+1.3) node[above=7pt,midway] {$\Delta T = m\Delta t$};

%
%
%
%
%
%
%
%
%
%
%
	\end{tikzpicture}
	\caption{\label{fig:time_grid} Uniformly spaced fine and coarse temporal meshes. $C$-points are present on both grids while $F$-points are only on the fine grid.}
	\end{center}
\end{figure}

\begin{figure}
\begin{center}
	\begin{tikzpicture}
\draw[line width=1.2pt, -, >=latex'](0,0) -- coordinate (x axis) (3,0) node[right] {}; 
\foreach \x in {0,1,2,3,4,5,6} \draw[ultra thick, black] (0.5*\x,0.125) -- (0.5*\x,-0.125) node[below] {$t_{\x}$} ;
\foreach \x in {0,3,6} \draw[ultra thick, black] (0.5*\x,0.225) -- (0.5*\x,-0.225) node[below] {} ;

\node (a) at (0.25,0.6) {$\boldsymbol{\Phi}_1$};
\node (a) at (0.75,0.6) {$\boldsymbol{\Phi}_2$};
\node (a) at (1.75,0.6) {$\boldsymbol{\Phi}_4$};
\node (a) at (2.25,0.6) {$\boldsymbol{\Phi}_5$};
\node (point_00) at (-0.1,0.2) {};
\node (point_01) at (0.6,0.2) {};
\node (point_10) at (0.4,0.2) {};
\node (point_11) at (1.1,0.2) {};

\node (point_20) at (1.4,0.2) {};
\node (point_21) at (2.1,0.2) {};
\node (point_30) at (1.9,0.2) {};
\node (point_31) at (2.6,0.2) {};

\draw[-{>[scale=0.7]}, ultra thick] (point_00) to[bend left = 40] (point_01);
\draw[-{>[scale=0.7]}, ultra thick] (point_10) to[bend left = 40] (point_11);

\draw[-{>[scale=0.7]}, ultra thick] (point_20) to[bend left = 40] (point_21);
\draw[-{>[scale=0.7]}, ultra thick] (point_30) to[bend left = 40] (point_31);

\draw[line width=1.2pt, -, >=latex'](4,0) -- coordinate (x axis) (7,0) node[right] {}; 
\foreach \x in {0,1,2,3,4,5,6} \draw[ultra thick, black] (4+0.5*\x,0.125) -- (4+0.5*\x,-0.125) node[below] {$t_{\x}$} ;
\foreach \x in {0,3,6} \draw[ultra thick, black] (4+0.5*\x,0.225) -- (4+0.5*\x,-0.225) node[below] {} ;

\node (a) at (5.25,0.6) {$\boldsymbol{\Phi}_3$};
\node (a) at (6.75,0.6) {$\boldsymbol{\Phi}_6$};

\node (point_00) at (4+0.9,0.2) {};
\node (point_01) at (4+1.6,0.2) {};

\node (point_10) at (4+2.4,0.2) {};
\node (point_11) at (4+3.1,0.2) {};

\draw[-{>[scale=0.7]}, ultra thick] (point_00) to[bend left = 40] (point_01);
\draw[-{>[scale=0.7]}, ultra thick] (point_10) to[bend left = 40] (point_11);

\node (a) at (1.5,-0.8) {$F$-relaxation};
\node (a) at (5.5,-0.8) {$C$-relaxation};

	\end{tikzpicture}
	\caption{\label{fig:relaxation} $F$- and $C$-relaxation for a grid with seven points and a temporal coarsening factor of three.}
	\end{center}
\end{figure}
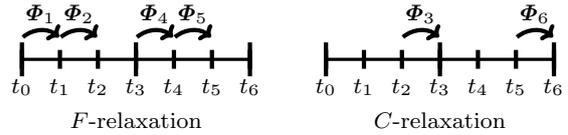

The MGRIT algorithm \cite{RDFalgout_etal_2014} is an iterative method that applies multigrid reduction (MGR) techniques \cite{Ries_MGR_methods} for solving time-dependent problems like \eqref{eq:problem_system}. The key practical aspect of MGRIT is its non-intrusiveness, i.\,e., the algorithm allows reusing existing time propagators and their integration into a parallel framework. For the construction of a multilevel algorithm, reusing the time integrators $\boldsymbol{\Phi}_i$, we require a definition of a coarse-grid system, a restriction and a prolongation operator between temporal grids and a relaxation scheme. Therefore, given a positive integer coarsening factor $m>1$, we define a splitting of the fine-grid points into $F$- and $C$-points, whereby every $m$-th point is a $C$-point and the others are $F$-points, resulting in a coarse time grid $T_{i_c} = i_c\Delta T, \; i=0,1,\dots,N_T$, with $N_T = N_t/m$ and time step $\Delta T = m \Delta t$ as depicted in Figure \ref{fig:time_grid}. 
For the coarse grid, we define a coarse-grid time integrator $\boldsymbol{\Phi}_{i_c}$ resulting from a re-discretization with time step $\Delta T$ of the problem; other choices, such as coarsening in the order of the discretization~\cite{Nielsen2018,Falgout_etal_jcomp2019} are also possible, but are not considered in this paper. We define two types of block relaxation schemes. The so-called $F$-relaxation performs a relaxation on $F$-points and propagates the solution from one $C$-point to all $F$-points up to the next $C$-point. Similarly, the $C$-relaxation propagates the solution \revision{from the preceding $F$-point} to a $C$-point. Both relaxation schemes, depicted in Figure \ref{fig:relaxation}, are highly parallel and can be combined to new relaxation schemes. For example, the $FCF$-relaxation, an $F$-relaxation followed by a $C$- and another $F$-relaxation, is typically used in the MGRIT algorithm. We define two grid transfer operators, injection as restriction and the ``ideal'' prolongation, for the transfer between temporal grid\revision{s}. \revision{The injection operator takes the values at the corresponding $C$-points in order to obtain values on the coarse grid.} The ideal interpolation is defined as \revision{the transpose of the} injection at $C$-points followed by an $F$-relaxation. \revision{That is, to obtain the values on the fine mesh the values at the $C$-points are obtained from the corresponding points on the coarse grid and these values are then propagated to the $F$-points in the current interval. The attribute ``ideal'' originates from the theory of algebraic multigrid as it is presented in \cite{Stueben2001}, as for this choice the Galerkin coarse grid operator is the Schur complement after reordering the system matrix according to the $C/F$-splitting.}

Using these components and the full approximation storage (FAS) framework \cite{Brandt_1977_FAS}, the two-level MGRIT-FAS algorithm \cite{Falgout_2014_FAS}, as extended by spatial coarsening in \cite{MR3716560}, can be written as Algorithm \ref{alg:mgrit}. \revision{Note that all statements of Algorithm \ref{alg:mgrit}, except for the solving in step 6, are highly parallel and can be executed simultaneously for several time steps}

\begin{algorithm}
\caption{\label{alg:mgrit}FAS-MGRIT($\mathcal{A}, \vec{u}, \mathbf{g}$)}
\begin{algorithmic}[1]
\State Apply $F$-relaxation to $\mathcal{A}_1(\textbf{u}^{(1)})= \textbf{g}^{(1)}$
\State For $0$ to $\gamma$:
\State \qquad Apply $CF$-relaxation to $\mathcal{A}_1(\textbf{u}^{(1)})= \textbf{g}^{(1)}$
\State Inject the approximation and its residual to the coarse grid: 
\Statex \qquad $\textbf{u}^{(2)} = \vec{R}_I(\textbf{u}^{(1)}), $
\Statex \qquad $\textbf{g}^{(2)} = \vec{R}_{I}(\textbf{g}^{(1)}-\mathcal{A}_1\textbf{u}^{(1)})$
\State If spatial coarsening:
\Statex \qquad $\textbf{u}^{(2)} = \vec{R}_s(\textbf{u}^{(2)})$ 
\Statex \qquad $\textbf{g}^{(2)} = \vec{R}_s(\textbf{g}^{(2)})$ 
\State Solve $\mathcal{A}_{2}(\textbf{v}^{(2)})= \mathcal{A}_{2}(\textbf{u}^{(2)}) + \textbf{g}^{(2)}$
\State Compute the error approximation: $\mathbf{e} = \mathbf{v}^{(2)} - \vec{u}^{(2)}$
\State If spatial coarsening:
\Statex \quad $\mathbf{e} = \vec{P}_s(e)$
\State Correct using ideal interpolation: $\textbf{u}^{(1)} = \textbf{u}^{(1)} + \vec{P}(\mathbf{e}) $
\end{algorithmic}
\end{algorithm}

In Algorithm \ref{alg:mgrit}, $\mathcal{A}_l(\textbf{u}^{(l)})= \textbf{g}^{(l)}$ denotes the space-time system of equations on levels $l=1,2$\revision{, which can be either linear or nonlinear,} with each block row corresponding to one time step. The relaxation scheme of the algorithm can be controlled by the parameter $\gamma$, and is, in general, for the MGRIT algorithm chosen as $\gamma = 1$, i.\,e., an $FCF$-relaxation scheme. But other choices are also possible. Note that the two-level variant with $F$-relaxation is equivalent to the parareal method \cite{JLLions_etal_2001a}. We denote the temporal grid transfer operators by $\mathbf{R}_I$ and $\mathbf{P}$, whereby $\mathbf{P}$ is the ideal prolongation and $\mathbf{R}_I$ the restriction operator based on injection. Additionally, we define $\vec{R}_s$ as spatial restriction and $\mathbf{P}_s$ as spatial prolongation (see Section \ref{sec:mgrit_spatial_coarsening}).

\begin{figure}
\begin{center}
\begin{tikzpicture}
\coordinate[label=left:] (c_0) at (0,0);
\coordinate[label=right:] (r_0_1) at (0.3,0.5);
\coordinate[label=right:] (r_0_2) at (0.6,1);
\coordinate[label=right:] (r_0_3) at (0.9,1.5);
\coordinate[label=right:] (r_0_4) at (1.2,2);

\coordinate[label=right:] (l_0_1) at (-0.3,0.5);
\coordinate[label=right:] (l_0_2) at (-0.6,1);
\coordinate[label=right:] (l_0_3) at (-0.9,1.5);
\coordinate[label=right:] (l_0_4) at (-1.2,2);

\draw[-, line width=2pt]  (c_0) -- (r_0_1);
\draw[-, line width=2pt] (r_0_1) -- (r_0_2);
\draw[-, line width=2pt] (r_0_2) -- (r_0_3);

\draw[-, line width=2pt]  (l_0_1) -- (c_0);
\draw[-, line width=2pt] (l_0_2) -- (l_0_1);
\draw[-, line width=2pt] (l_0_3) -- (l_0_2);

\fill (c_0) circle (3pt);
\fill (r_0_1) circle (3pt);
\fill (r_0_2) circle (3pt);
\fill (r_0_3) circle (3pt);
\fill (l_0_1) circle (3pt);
\fill (l_0_2) circle (3pt);
\fill (l_0_3) circle (3pt);

\coordinate[label=left:] (c_1) at (3,0);
\coordinate[label=right:] (r_1_1) at (3.3,0.5);
\coordinate[label=right:] (r_1_2) at (3.6,0);
\coordinate[label=right:] (r_1_3) at (3.9,0.5);
\coordinate[label=right:] (r_1_4) at (4.2,1);
\coordinate[label=right:] (r_1_5) at (4.5,0.5);
\coordinate[label=right:] (r_1_6) at (4.8,0);
\coordinate[label=right:] (r_1_7) at (5.1,0.5);
\coordinate[label=right:] (r_1_8) at (5.4,1);
\coordinate[label=right:] (r_1_9) at (5.7,1.5);
\coordinate[label=right:] (r_1_10) at (6,2);

\coordinate[label=right:] (l_1_1) at (2.7,0.5);
\coordinate[label=right:] (l_1_2) at (2.4,1);
\coordinate[label=right:] (l_1_3) at (2.1,1.5);
\coordinate[label=right:] (l_1_4) at (1.8,2);

\draw[-, line width=2pt]  (c_1) -- (r_1_1);
\draw[-, line width=2pt]  (r_1_1) -- (r_1_2);
\draw[-, line width=2pt]  (r_1_2) -- (r_1_3);
\draw[-, line width=2pt]  (r_1_3) -- (r_1_4);
\draw[-, line width=2pt]  (r_1_4) -- (r_1_5);
\draw[-, line width=2pt]  (r_1_5) -- (r_1_6);
\draw[-, line width=2pt]  (r_1_6) -- (r_1_7);
\draw[-, line width=2pt]  (r_1_7) -- (r_1_8);
\draw[-, line width=2pt]  (r_1_8) -- (r_1_9);

\draw[-, line width=2pt]  (l_1_1) -- (c_1);
\draw[-, line width=2pt] (l_1_2) -- (l_1_1);
\draw[-, line width=2pt] (l_1_3) -- (l_1_2);

\fill (c_1) circle (3pt);
\fill (r_1_1) circle (3pt);
\fill (r_1_2) circle (3pt);
\fill (r_1_3) circle (3pt);
\fill (r_1_4) circle (3pt);
\fill (r_1_5) circle (3pt);
\fill (r_1_6) circle (3pt);
\fill (r_1_7) circle (3pt);
\fill (r_1_8) circle (3pt);
\fill (r_1_9) circle (3pt);
\fill (l_1_1) circle (3pt);
\fill (l_1_2) circle (3pt);
\fill (l_1_3) circle (3pt);
\end{tikzpicture}
	\caption{\label{fig:v_f_cycle} Structure of $V$- and $F$-cycles for four grid levels.}
	\end{center}
\end{figure}
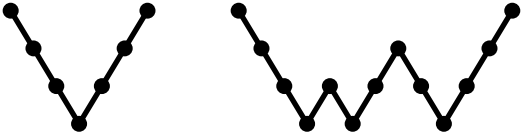

The two-level MGRIT-FAS algorithm can be extended to a multilevel setting by applying the two-level method recursively to the system in step $6$. Depending on the number of recursive calls, different multilevel schemes, e.\,g., $V$-, $F$- or $W$-cycles as depicted in Figure \ref{fig:v_f_cycle}, can be defined. Note that, in exact arithmetic the MGRIT algorithm with $F$- or $FCF$-relaxation solves for the discrete fine-grid solution in $N_t/m$ or $N_t/\left(2m\right)$ iterations \cite{RDFalgout_etal_2014}, respectively.

In Algorithm \ref{alg:mgrit}, an initial guess $\vec{u}$ is needed. In general, if prior knowledge of the solution $\vec{u}$ exists, this knowledge should be used as an initial guess. If no such prior knowledge exist, an improved initial guess can be computed by using the nested iterations \cite{GDahlquist_nes_it,LKronsjo_nes_it} strategy. The idea of nested iterations is to compute an initial approximation on the coarsest grid and interpolate the approximation to the finer grids, applying one $V$-cycle per level. This procedure results in the same structure as the $F$-cycle, with skipping the down-cycle at the beginning.

\subsection{MGRIT with spatial coarsening}
\label{sec:mgrit_spatial_coarsening}

Spatial coarsening is added to the MGRIT-FAS algorithm by the spatial restriction operator $\vec{R}_s$ in line $5$ and the spatial prolongation operator $\vec{P}_s$ in line $8$. Note, that spatial coarsening is an extra option that can be applied after semi-coarsening in time and, therefore, does not touch the non-intrusive character of the algorithm. So both operators have to be specified manually, but are independent of the time-integration scheme. Obviously, the spatial grid transfer operators directly influence the convergence behavior of the algorithm and should be chosen wisely.

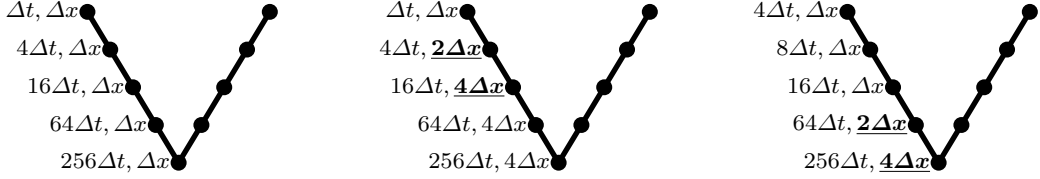
\begin{figure*}
\begin{center}
\begin{tikzpicture}
\coordinate[label=left:{$256\Delta t, \Delta x$}] (c_0) at (0,0);
\coordinate[label=right:] (r_0_1) at (0.3,0.5);
\coordinate[label=right:] (r_0_2) at (0.6,1);
\coordinate[label=right:] (r_0_3) at (0.9,1.5);
\coordinate[label=right:] (r_0_4) at (1.2,2);

\coordinate[label=left:{$64\Delta t, \Delta x$}] (l_0_1) at (-0.3,0.5);
\coordinate[label=left:{$16\Delta t, \Delta x$}] (l_0_2) at (-0.6,1);
\coordinate[label=left:{$4\Delta t, \Delta x$}] (l_0_3) at (-0.9,1.5);
\coordinate[label=left:{$\Delta t, \Delta x$}] (l_0_4) at (-1.2,2);

\draw[-, line width=2pt]  (c_0) -- (r_0_1);
\draw[-, line width=2pt]  (r_0_1) -- (r_0_2);
\draw[-, line width=2pt] (r_0_2) -- (r_0_3);
\draw[-, line width=2pt] (r_0_3) -- (r_0_4);

\draw[-, line width=2pt]  (l_0_1) -- (c_0);
\draw[-, line width=2pt] (l_0_2) -- (l_0_1);
\draw[-, line width=2pt] (l_0_3) -- (l_0_2);
\draw[-, line width=2pt] (l_0_4) -- (l_0_3);

\fill (c_0) circle (3pt);
\fill (r_0_1) circle (3pt);
\fill (r_0_2) circle (3pt);
\fill (r_0_3) circle (3pt);
\fill (r_0_4) circle (3pt);
\fill (l_0_1) circle (3pt);
\fill (l_0_2) circle (3pt);
\fill (l_0_3) circle (3pt);
\fill (l_0_4) circle (3pt);

\coordinate[label=left:{$256\Delta t, 4\Delta x$}] (c_1) at (5,0);
\coordinate[label=right:] (r_1_1) at (5.3,0.5);
\coordinate[label=right:] (r_1_2) at (5.6,1);
\coordinate[label=right:] (r_1_3) at (5.9,1.5);
\coordinate[label=right:] (r_1_4) at (6.2,2);

\coordinate[label=left:{$64\Delta t, 4\Delta x$}] (l_1_1) at (4.7,0.5);
\coordinate[label=left:{$16\Delta t, \underline{\boldsymbol{4\Delta x}}$}] (l_1_2) at (4.4,1);
\coordinate[label=left:{$4\Delta t, \underline{\boldsymbol{2\Delta x}}$}] (l_1_3) at (4.1,1.5);
\coordinate[label=left:{$\Delta t, \Delta x$}] (l_1_4) at (3.8,2);

\draw[-, line width=2pt]  (c_1) -- (r_1_1);
\draw[-, line width=2pt]  (r_1_1) -- (r_1_2);
\draw[-, line width=2pt] (r_1_2) -- (r_1_3);
\draw[-, line width=2pt] (r_1_3) -- (r_1_4);

\draw[-, line width=2pt]  (l_1_1) -- (c_1);
\draw[-, line width=2pt] (l_1_2) -- (l_1_1);
\draw[-, line width=2pt] (l_1_3) -- (l_1_2);
\draw[-, line width=2pt] (l_1_4) -- (l_1_3);

\fill (c_1) circle (3pt);
\fill (r_1_1) circle (3pt);
\fill (r_1_2) circle (3pt);
\fill (r_1_3) circle (3pt);
\fill (r_1_4) circle (3pt);
\fill (l_1_1) circle (3pt);
\fill (l_1_2) circle (3pt);
\fill (l_1_3) circle (3pt);
\fill (l_1_4) circle (3pt);

\coordinate[label=left:{$256\Delta t, \underline{\boldsymbol{4\Delta x}}$}] (c_2) at (10,0);
\coordinate[label=right:] (r_2_1) at (10.3,0.5);
\coordinate[label=right:] (r_2_2) at (10.6,1);
\coordinate[label=right:] (r_2_3) at (10.9,1.5);
\coordinate[label=right:] (r_2_4) at (11.2,2);

\coordinate[label=left:{$64\Delta t, \underline{\boldsymbol{2\Delta x}}$}] (l_2_1) at (9.7,0.5);
\coordinate[label=left:{$16\Delta t, \Delta x$}] (l_2_2) at (9.4,1);
\coordinate[label=left:{$8\Delta t, \Delta x$}] (l_2_3) at (9.1,1.5);
\coordinate[label=left:{$4\Delta t, \Delta x$}] (l_2_4) at (8.8,2);

\draw[-, line width=2pt]  (c_2) -- (r_2_1);
\draw[-, line width=2pt]  (r_2_1) -- (r_2_2);
\draw[-, line width=2pt] (r_2_2) -- (r_2_3);
\draw[-, line width=2pt] (r_2_3) -- (r_2_4);

\draw[-, line width=2pt]  (l_2_1) -- (c_2);
\draw[-, line width=2pt] (l_2_2) -- (l_2_1);
\draw[-, line width=2pt] (l_2_3) -- (l_2_2);
\draw[-, line width=2pt] (l_2_4) -- (l_2_3);

\fill (c_2) circle (3pt);
\fill (r_2_1) circle (3pt);
\fill (r_2_2) circle (3pt);
\fill (r_2_3) circle (3pt);
\fill (r_2_4) circle (3pt);
\fill (l_2_1) circle (3pt);
\fill (l_2_2) circle (3pt);
\fill (l_2_3) circle (3pt);
\fill (l_2_4) circle (3pt);
\end{tikzpicture}
	\caption{\label{fig:coarsening_strategies} Five-level MGRIT $V$-cycle algorithm with three spatial grids and different space-time coarsening strategies. The left figure shows the MGRIT algorithm with coarsening only in time, the middle one illustrates direct spatial coarsening, and at right, the delayed spatial coarsening strategy is applied.}
	\end{center}
\end{figure*}

In a multilevel setting, multiple coarsening strategies can be chosen depending on the number of available grids in time and space. Both coarsening dimensions are independent of each other and, therefore, can be applied in different ways. Assuming that more temporal grids than spatial grids are available, the most cost efficient approach is to use a ``direct'' spatial coarsening strategy, starting with spatial coarsening as early as possible. This strategy would \revision{reduce the computational costs per iteration for all spatial coarse level problems}, but studies \cite{MR3716560} have \revision{demonstrated} that this ``direct'' strategy can degrade the convergence behavior of the MGRIT algorithm for linear and nonlinear problems. To overcome the degradation in MGRIT convergence, in \cite{MR3716560} a ``delayed'' spatial coarsening strategy is proposed. This strategy applies spatial coarsening as late as possible, i.\,e., only on the coarsest time grids. Figure \ref{fig:coarsening_strategies} illustrates the different coarsening strategies for a five-level $V$-cycle for the case of three spatial grids, characterized by spatial grid spacings of $\Delta x$, $2\Delta x$, and $4\Delta x$, i.\,e., assuming factor-two coarsening in space. The temporal coarsening factor is chosen to be $m=4$, i.\,e., the time step on each temporal grid is given by $4^{l-1}\Delta t, ~l > 0$. The left $V$-cycle represents the MGRIT algorithm with coarsening only in the time dimension. In the middle, the direct spatial coarsening strategy is shown, where spatial coarsening  is applied on the first and second time levels. Finally, the right $V$-cycle illustrates the delayed strategy, with spatial coarsening starting on the third level to make use of the coarsest space grid on the coarsest level. 

\section{Simulation of induction machines}
\label{sec:sim_induction_machine}

The simulation of electromagnetic fields is based on the solution of Maxwell's equations \cite{jackson_classical_1999}, but, for a variety of applications, it is sufficient to solve a (quasistatic) approximation. In this section, we introduce one of these simplifications, the so-called eddy current problem, that is typically used in the simulation of energy converters. Further, we introduce a numerical model of an induction machine that is used for computations in Sections~\ref{sec:linear_model} and \ref{sec:nonlinear_model}.

\subsection{Induction machine model}
\label{sec:induction_machine_model}
\revision{
The standard approach in the simulation of electrical machines is to neglect the displacement current density with respect to the source current density. The resulting magnetoquasistatic approximation of Maxwell's equations is the well-known eddy current problem.

Let $\Omega$ denote a 2D cross-section in $x,y$-coordinates of the 3D geometry (with depth $\ell_z$). Then, the problem can be solved in terms of the $z$-component of the magnetic vector potential ${A}_z:\Omega \times \timeint \rightarrow \mathbb{R}$. We split the spatial domain into stator and rotor $\bar\Omega=\bar\Omega_\text{1}\cup\bar\Omega_\text{2}$ and introduce separate unknowns ${A}_k={A}_z|_{\Omega_k}$ on each domain ($k=1,2$). The problem reads for the stator
\begin{subequations}\label{eq:system_eddy_current}
\begin{align}
		- \nabla \cdot \left(\nu \nabla {A}_1\right) -\sum_{s=1}^{n_{1}} \vec{\chi}_{s} i_{1,s} & = 0
		\label{eq:system_eddy_current_a} 
		\\
		\frac{d}{dt}\int_{\Omega} {\chi}_{s}{A}_1\; \ell_z\; dS  + R_{s}(\vec{i}_{1}) & = v_{1,s},
		\label{eq:system_eddy_current_b}
\intertext{where $s=1,\ldots, n_{1}$, and for the rotor}
		\sigma \partial_t {A}_2 - \nabla \cdot \left(\nu \nabla{A}_2\right) -\sum_{s=1}^{n_{2}} \sigma\vec{\zeta}_{s} v_{2,s} & = 0
		\label{eq:system_eddy_current_c} 
		\\
		-\frac{d}{dt}\int_{\Omega} \vec{\zeta}_{s}\sigma{A}_2\; \ell_z\; dS  + G_s(\vec{v}_{2},\vec{i}_{2}) & = 0
		\label{eq:system_eddy_current_d}
		\\
		{L}_s(\vec{v}_{2},\frac{d}{dt}\vec{i}_{2}) & =0,
\end{align}
where $s=1,\ldots,n_{2}$, with interface conditions on $\Omega_1\cap\Omega_2$ within the air-gap for all $t\in\timeint$,
\begin{align}
		{A}_1(\vec{x}_1,t) &= {A}_2(\vec{x}_2,t)\\
		\vec{n}\cdot\nabla{A}_1(\vec{x}_1,t) &= \vec{n}\cdot\nabla{A}_2(\vec{x}_2,t),
\end{align}
\end{subequations}
where $\vec{x}_2=\vec{r}(\vec{x}_1,\theta)$ imposes the rotation by angle $\theta$. The problem is completed by the Dirichlet boundary condition $A_z = 0$ on $\partial\Omega$ and the initial value ${A}_z\left(\vec{x},t_0\right) = {A}_0(\vec{x})$. 
The electrical conductivity $\sigma(\vec{x}) \geq 0$ and the (isotropic, nonlinear) magnetic reluctivity $\nu(\vec{x},|\nabla{A}|) > 0$ encode the geometry. The current and voltage distribution functions ${\chi}_{s},{\zeta}_{s}: \Omega_s \rightarrow \mathbb{R}$ model the stator coils by stranded conductors and the rotor bars by interconnected solid conductors \cite{schoeps2013windingfunctions} in terms of resistive circuits $R_{s}$ and $G_{s}$, inductive circuits $L_s$, voltages $\vec{v}_k=[v_{k,1},\ldots,v_{k,n_k}]^\top$, and currents $\vec{i}_k=[i_{k,1},\ldots,i_{k,n_k}]^\top$ for $k=1,2$.
}

A common quantity of interest in the design of electrical machines are the Joule losses
\begin{equation}
	\label{eq:losses}
	P_\textrm{loss}=\int_\Omega \partial_t\revision{A_z}\;\sigma\partial_t\revision{A_z\;\ell_z}\;  d\Omega.
\end{equation}

To include the rotation of the rotor, the problem above is extended by an additional motion equation. Movements are determined by the mechanical equation 
\begin{subequations}\label{eq:system_motion}
	\begin{alignat}{2}
	\omega(t)&= \frac{d\theta(t)}{dt}, \quad &&t \in \timeint, \label{eq:system_motion_a}\\
		I\frac{d^2\theta(t)}{dt^2}+C\frac{d\theta(t)}{dt} &= T_{\text{mag}}\left(\revision{A_z}\right) \quad &&\text{in} \; \timeint, \label{eq:system_motion_b} \\ 
		\theta\left(t_0\right)&=\theta_0, \label{eq:system_motion_c} \\
		\omega(t_0) &= \omega_0, \label{eq:system_motion_d}
	\end{alignat}
\end{subequations}
where $\omega$ is the angular velocity of the rotor, $I$ denotes the moment of inertia, and $C$ is the friction coefficient, respectively. The torque $T_\text{mag}$ defines the mechanical excitation of system \eqref{eq:system_motion} given by the magnetic field. \revision{Note that there are several common special cases, i.\,e., a fixed rotor ($\theta = \text{const}$) without movement and given rotor speed (synchronous or asynchronous).}

Discretizing \eqref{eq:system_eddy_current} in space using finite elements with $n_a$ degrees of freedom \revision{and using the moving band approach \cite{FerreiraDaLuz_etal_moving_band} for motion} yields 
a system of equations of the form
\begin{subequations}\label{eq:semi_discrete_system}
	\begin{align}
		M \mathbf{u}'(t) + K(\mathbf{u}(t)) \mathbf{u}(t) &= \mathbf{f}(t), \quad t \in \timeint \label{eq:semi_discrete_system_a}\\
		\mathbf{u}(t_0) &=\mathbf{u}_0, \label{eq:semi_discrete_system_b}
	\end{align}
\end{subequations}
\revision{with initial condition $\mathbf{u}_0 \in \mathbb{R}^n$ and unknown $\mathbf{u}^\top=[\mathbf{a}^\top, \mathbf{i}_1^\top, \mathbf{i}_2^\top, \mathbf{v}_2^\top, \theta, \omega] : \timeint \rightarrow \mathbb{R}^n$. The solution $\mathbf{u}$ consists of the magnetic vector potential $\vec{a}$, currents $\vec{i}_1$, $\vec{i}_2$, voltages $\vec{v}_2$, the rotor angle $\theta$, and the rotor's angular velocity $\omega$. The stator voltages $v_{1,s}(t)$ drive the problem and, thus, determine the right-hand side $\mathbf{f}(t)$.} Integrating the semi-discrete system \eqref{eq:semi_discrete_system} using the backward Euler method results in a system of the form
\begin{subequations}\label{eq:discrete_system}
	\begin{align}
		\left( \frac{1}{\Delta t}M + K\left( \mathbf{u_{i}} \right) \right) \mathbf{u_{i}} &= \mathbf{f_{i}}+ \frac{1}{\Delta T}M \mathbf{u}_{i-1}, \label{eq:discrete_system_a}\\
		\vec{u}_0 &=\vec{u}(t_0). \label{eq:discrete_system_b}
	\end{align}
\end{subequations}
This field/circuit coupled problem consists of dif\-ferential-algebraic equations (DAEs) of index-1 \cite{Bartel_2011aa,Cortes-Garcia_2019aa} depending on the circuit topology. A discussion of the index-1 DAE aspects in a parallel-in-time setting can be found in \cite{Schops_2018aa}\revision{, i.e., if the mass matrix $M$ is constant and the backward Euler method is used then there is no need for explicitly enforcing consistency of the (algebraic) initial values}.

\subsection{Pulse-width-modulation}
\label{sec:pwm}

\begin{figure}
	\begin{center}
		\includegraphics[width=0.4\textwidth]{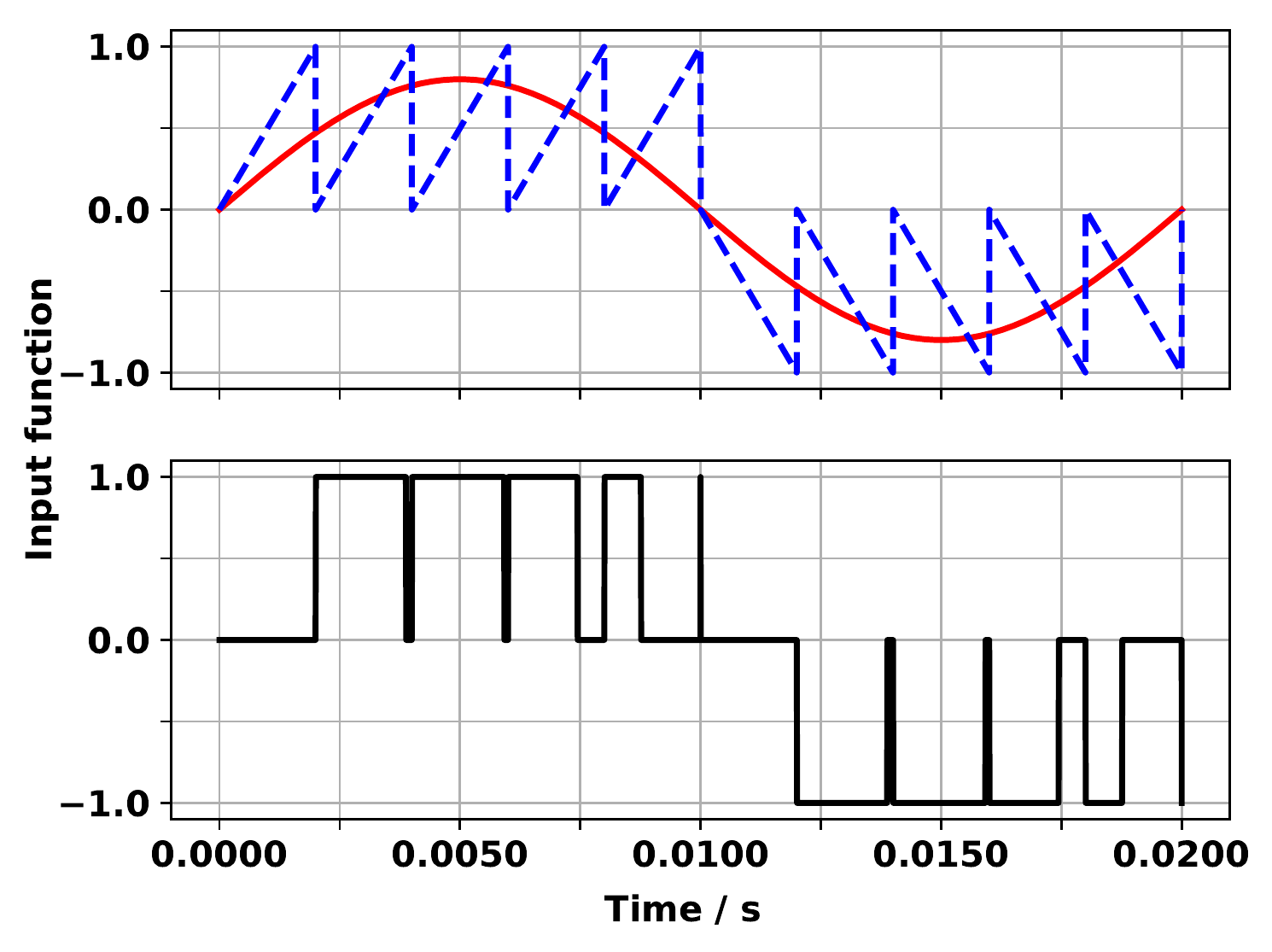}
			\caption{Illustration of a PWM signal (black), generated by a sine wave reference signal (red) and a sawtooth carrier signal (blue).}
		\label{fig:pwm}
	\end{center}
\end{figure}

Electrical machines are commonly driven by PWM excitations since they have technical advantages for the power electronics that is implemented to control the machine. While multiple forms of PWM exists \cite{black1953modulation}, in this paper, we consider a PWM signal that is produced by comparing a reference wave with a triangular wave. Idealizing the circuitry to generate the PWM, i.\,e. diodes are replaced by switches,  results in the following excitation
\begin{equation}
	b\left(t\right) = \text{sgn} \left[ r \left( t \right) - c \left( t \right) \right ], \quad t\in\mathbb{R},
\end{equation}
with reference signal $r\left(t\right)$ and where $c\left(t\right)$ denotes the carrier signal. Furthermore, the modulation factor is defined as the ratio of the amplitude of the reference signal and the amplitude of the carrier signal. Figure \ref{fig:pwm} shows an example of a PWM signal. The red line illustrates a sine reference signal of $50$Hz with an amplitude of $0.8$ and the carrier signal, modeled by a sawtooth carrier \cite{Sun2012_pwm} of $500$Hz and with an amplitude of one, is shown as a blue line. The resulting PWM signal with switching frequency of $500$Hz is presented in the lower subplot of Figure \ref{fig:pwm}. 

\subsection{Numerical model}
\label{sec:numerical_model}

\begin{figure}
	\begin{center}
		\includegraphics[width=0.4\textwidth]{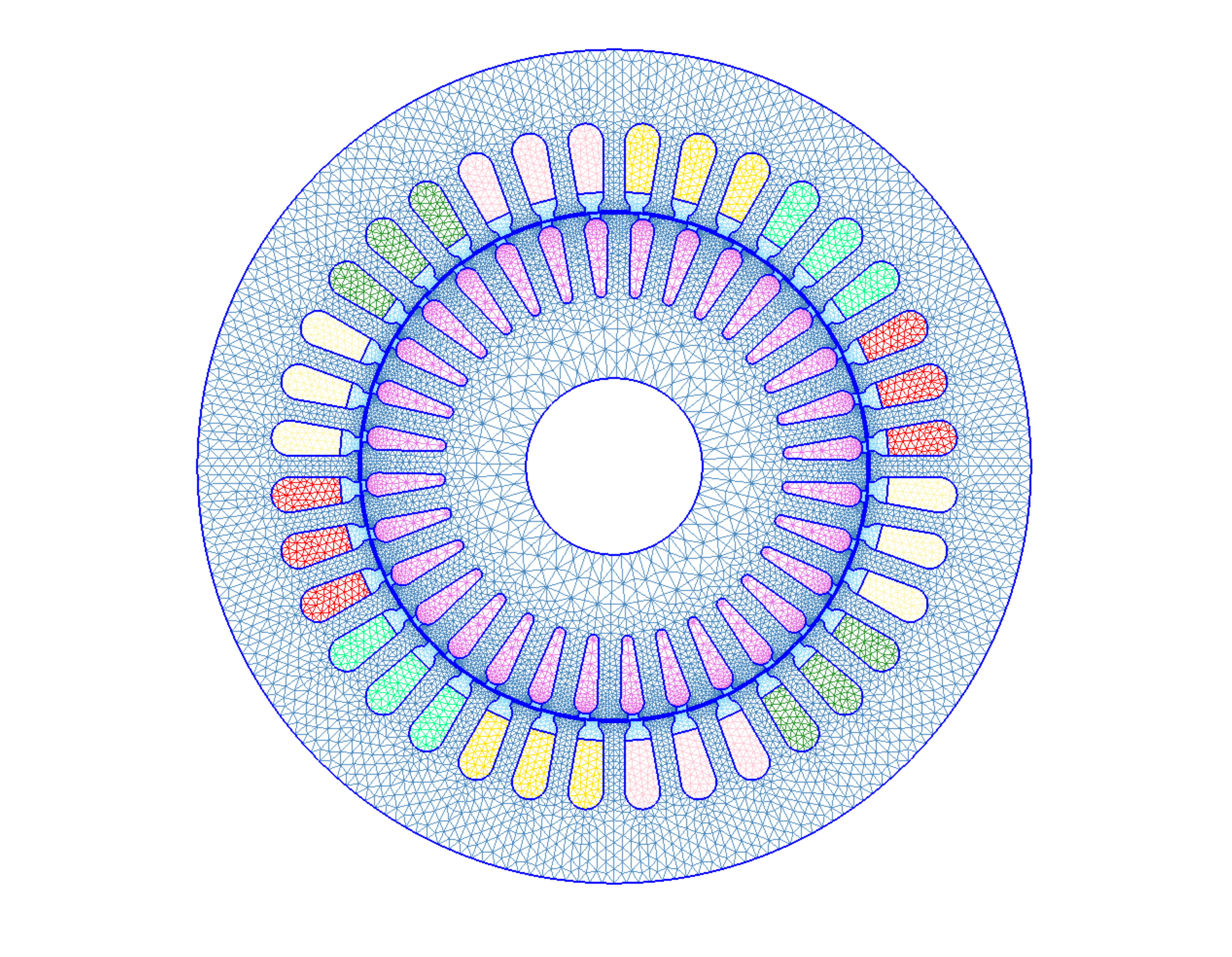}
			\caption{Mesh view of the four-pole squirrel cage machine model ``im\_3\_kw'' of \cite{JGyselinck2001_multi_slice}.}
		\label{fig:machine}
	\end{center}
\end{figure}

We model the semi-discrete eddy-current problem \eqref{eq:semi_discrete_system}, supplied by a three-phase PWM voltage source, using the multi-slice FE model ``im\_$3$\_kw'' of an electrical machine first presented in \cite{JGyselinck2001_multi_slice} and modified in \cite{Gander_2019aa} for the usage of a PWM voltage source. In more detail, the two-dimensional model ``im\_$3$\_kw'', depicted in Figure~\ref{fig:machine}, of a  four-pole $3$kW squirrel cage induction machine with no-load operation condition is used for numerical computations. Further, as typical for this kind of symmetric models, we consider only a quarter of the machine with periodic boundaries to reduce simulation costs. The machine is supplied by a three-phase PWM voltage source
\begin{equation} \label{eq:machine_pwm}
	v_s^p(t) = \text{sign} \left[ r_s \left( t \right) - c_p(t)\right], \qquad t \in \timeint
\end{equation}
with reference signals $r_s\left( t \right)$, $s=1,2,3$, and carrier signal $c_p\left( t \right)$ with $p$ pulses. The reference signal for the three phases is given by
\begin{equation} \label{eq:machine_reference_signal}
	r_s \left( t \right) = \sin \left( \frac{2\pi}{T}t+\varphi_s \right),
\end{equation}  
with $\varphi_1=0$, $\varphi_2=\-2/3\pi$, $\varphi=-4/3\pi$ and electric period $T$. The carrier signal, given by a bipolar trailing-edge modulation using a sawtooth carrier \cite{Sun2012_pwm}, is defined by
\begin{equation} \label{eq:machine_carrrier_signal}
 	c_p(t) = 2 \left( \frac{p}{T}t - \lfloor\frac{p}{T}t\rfloor \right) -1.
\end{equation}
For our simulation we consider a frequency of $50$ Hz, corresponding to an electric period of $T=0.02$ s, and choose $p=400$ pulses for one period, resulting in a PWM voltage source of $20$ kHz. The model provides a linear, i.\,e., $\nu$ does not depend on $\vec{A}$, and a nonlinear version of the induction machine. 
In the nonlinear setting, the voltage source $v_s^p$, defined in \eqref{eq:machine_pwm}, is multiplied with the time-dependent factor 
\begin{equation}
	\alpha\left(t\right)=0.5\left(1-\cos  \left( \pi \; t/2 \right) \right)
\end{equation}
in the first two periods, i.\,e., for $t \in \left[0,2T\right]$, which reduces the transient behavior of the motor and allows to obtain the steady-state more quickly \cite{JGyselinck2001_multi_slice}. 

The model computes the value of Joule losses $P_\textrm{loss}$ for each time point, i.\,e., the loss in the transfer of electrical energy. The size of Joule losses indicates a significant quality in the prototyping of designs, and, therefore, it is one possible stopping criterion for iterative methods.

\subsection{Spatial discretization}
\label{subsec:spatial_discretization}

\begin{figure}
	\centerline{\includegraphics[width=.32\linewidth]{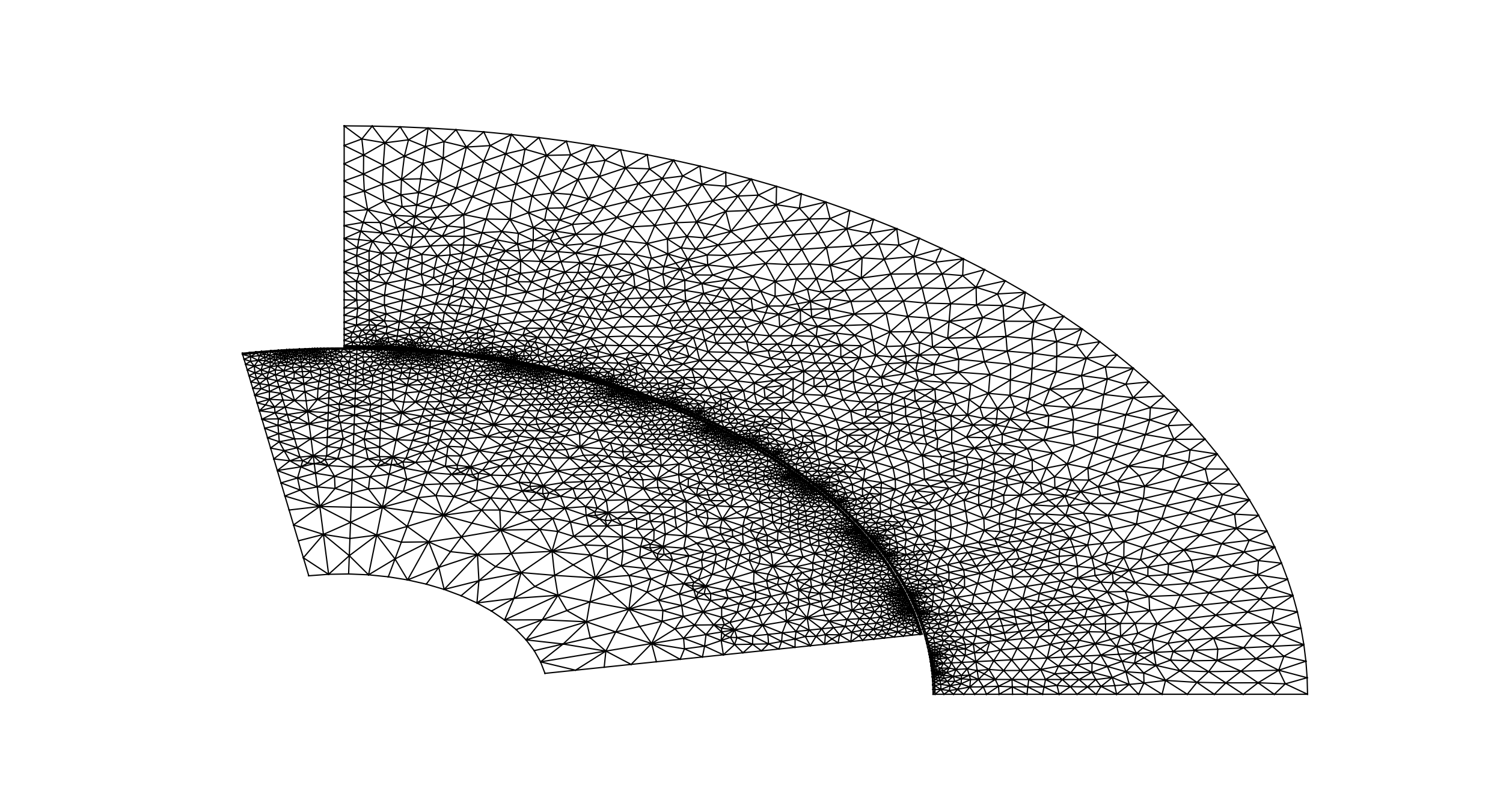}\quad
	\includegraphics[width=.32\linewidth]{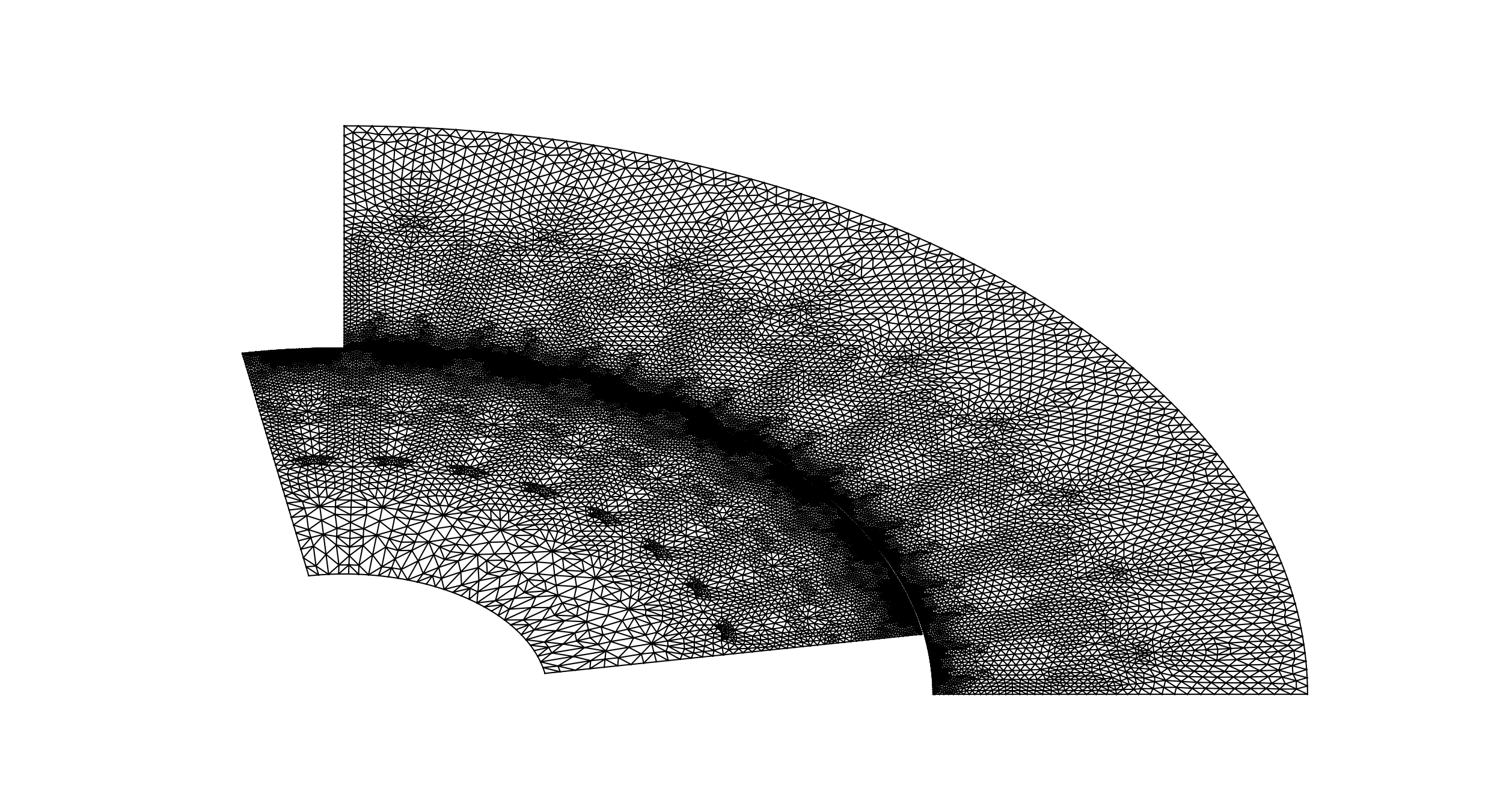}\quad
	\includegraphics[width=.32\linewidth]{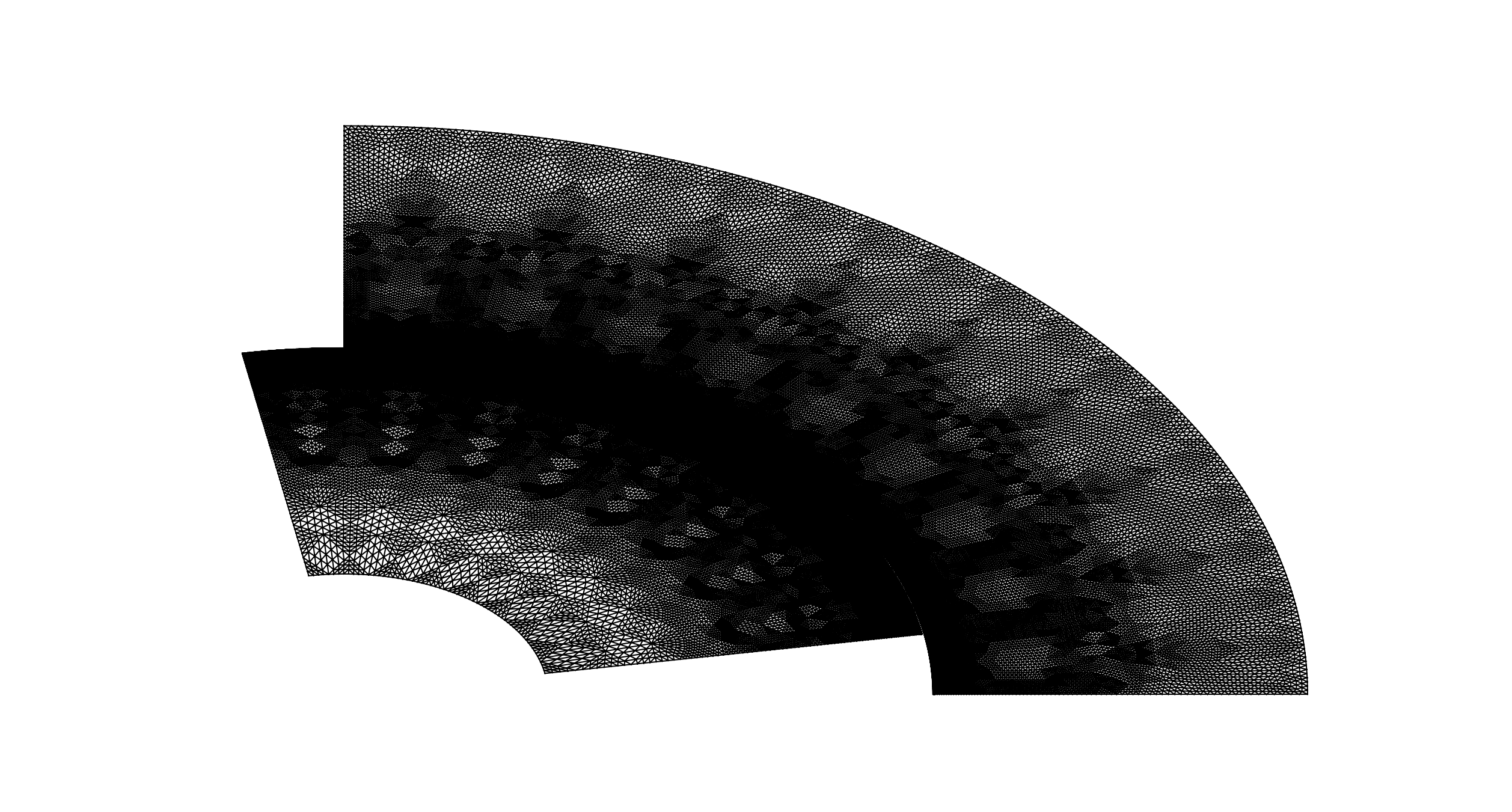}}
%
	\caption{Three spatial discretizations of one quarter of the machine generated by Gmsh. Finer meshes $\mesh_2$ (middle) and $\mesh_1$ (right) are created by refinements of the left mesh $\mesh_3$.}
	\label{fig:discretizations}
\end{figure}

Gmsh \cite{gmsh, gmsh_website} is used to generate mesh representations of the model. We consider a hierarchy of three nested spatial grids, illustrated in Figure \ref{fig:discretizations}. The coarsest model, defined on the mesh $\mesh_3$, has $n_a=4449$ degrees of freedom. Two more accurate discretizations are defined on grids $\mesh_2$  and $\mesh_1$ that are constructed by refinements of the coarsest mesh, resulting in grids with $n_a=17{,}496$ and $n_a=69{,}384$ degrees of freedom, respectively. 

We consider the standard finite element interpolation or nodal interpolation as the grid transfer operator between the spatial meshes. To avoid numerical instabilities at the boundaries between the stator and the rotor, we apply the standard finite element interpolation separately for both regions. Recall from Section \ref{sec:induction_machine_model} that in our model problem, only the magnetic vector potential is a function that depends on space, and, therefore, only the magnetic vector potential is discretized on the different spatial meshes, whereas currents, rotor angle, and the rotor's angular velocity are functions that are independent of space. Taking into account that System \eqref{eq:semi_discrete_system} consists of grid-dependent and grid-independent unknowns, we define the spatial interpolation and restriction operators as block operators with two blocks as follows: For the interpolation of grid-dependent unknowns, we use the standard finite element interpolation, whereas the grid-independent unknowns are injected.

\subsection{Implementation details}
\label{subsec:implementation_details}

The MGRIT algorithm is implemented in PyMGRIT \cite{PyMGRIT}, a Python package using the Message Passing Interface (MPI). The implementation conforms to the non-intrusive character of the MGRIT method and requires an implementation of the time stepper routine by the user. For the machine model, the time stepper routine calls the GetDP \cite{PDular_etal_1998,getdp, getdp_website} library for performing the simulation of the electrical machine which implements the implicit Euler method. As we aim at demonstrating the benefits of parallelization in time, for the implementation on a distributed memory computer, we consider a domain-decomposition approach only in time, i.\,e., only the time interval is distributed across processors. All tests were performed on an Intel Xeon Phi Cluster consisting of $272$ $1.4$ GHz Intel Xeon Phi processors.

\subsection{Storage requirements}
\label{subsec:storage_requirements}

The MGRIT algorithm requires storing solution values at all $C$-points. On coarse levels, two additional vectors are stored: one for a restricted copy of the fine-grid solution and one for the FAS right-hand-side. Denoting the number of MGRIT levels by $L$ and the storage cost of a vector on level $l$ by $s_l$, $l=0,\ldots,L-1$, the total storage requirement of MGRIT per temporal processor is given by~\cite{MR3716560}
\[
S_{\text{MGRIT}} \approx \sum_{l=0}^{L-1} \left\lceil \frac{N_t}{m^{l+1}p} \right\rceil s_l + \sum_{l=1}^{L-1} 2 \left\lceil \frac{N_t}{m^lp} \right\rceil s_l,
\]
where $p$ is the number of processors in time. The first term corresponds to the storage requirements of solution values, and the second term is that of the additional vectors on the coarse levels. Note that spatial coarsening reduces the size of all vectors on coarse grids.

In comparison, the sequential time-stepping algorithm requires storing only one solution value, corresponding to a total storage requirement of $s_0$. Thus, the total storage cost of the MGRIT algorithm is larger than that of sequential time-stepping by a factor of $S_{\text{MGRIT}}/s_0$. However, note that this factor decreases with increasing the number of processors and in the case of spatial coarsening. 

\section{Linear material model}
\label{sec:linear_model}

In this section, we investigate the behavior of the\linebreak MGRIT algorithm for the eddy current problem of a two-dimensional induction machine with linear material characteristics. 
Especially, we are interested in the behavior of adding spatial coarsening to the MGRIT algorithm for system \eqref{eq:semi_discrete_system}, containing scalar-valued as well as grid-dependent unknowns. For this study, we choose a linear model of the induction machine, also available in the model ``im\_$3$\_kw'', to reduce computational costs. In Section \ref{sec:nonlinear_model}, we transfer our results to the nonlinear model of the machine. 

\subsection{Algorithmic parameters}
\label{subsec:linear_algorithmic_parameters}

The MGRIT-FAS algorithm has a large number of parameters and variations, considering $V$-, $F$- and $W$-cycles, the number of temporal grid levels, coarsening strategies in space and time, number and type of relaxation schemes, and so forth. Based on results presented in the literature as well as on personal experiences, we consider a two- and a five-level MGRIT variant in this study. Further, we choose $V$- and $F$-cycles. $V$-cycles offer better parallel scaling than $F$-cycles, but $F$-cycles have better convergence behavior due to additional work on coarse temporal grids. We consider $FCF$-relaxation for all variants, and, additionally, $F$-relaxation is used in the two-level setting as well as in MGRIT $F$-cycles.

The main goal of the MGRIT algorithm is to reduce the runtime of the simulation by making use of parallel computer architectures. To achieve the best runtime results for our model problem, the temporal coarsening strategy of the algorithm has to be adapted to the number of available processors. Therefore, we choose a non-uniform temporal coarsening strategy directly connected to the number of points on the finest time grid and the number of processors. On the one hand, a large coarsening factor reduces the number of points for all following time grids and, thereby, the costs. On the other hand, to make use of all available processors, the coarsening factor cannot be chosen too large. As a compromise and based on experience with MGRIT for eddy current problems in a parallel setting \cite{Friedhoff_2020aa,friedhoff2019multigridreductionintime}, we choose a large coarsening factor on the first level, such that the number of points on the second level corresponds to the number of processors available. On coarse levels, a small coarsening factor is used. For example, consider a fine temporal grid with $129$ points, $16$ processors available and a five-level MGRIT method. The first temporal coarsening factor is chosen to be eight, such that the work of the first grid is distributed evenly among all processors; here, eight time points per processor. For the coarser levels, a coarsening factor of two is used to reduce the number of time points per grid, such that the coarsest grid has only three points. With this strategy, the first level is perfectly parallelized and the amount of work on coarser levels that cannot be parallelized is minimized. Note that this is a heuristic choice purely based on distributing computational work evenly among a given number of processors. Due to the high computational cost of the external solver and to the use of moderate numbers of processors, it is reasonable to neglect communication costs for this problem. Optimization of the temporal coarsening strategy that takes additional parameters that affect parallel performance, such as communication costs and MGRIT convergence behavior into account, is a topic of future work.

In MGRIT multilevel variants, we apply the direct and delayed spatial coarsening strategies described in Section \ref{sec:mgrit_spatial_coarsening}, with the interpolation and restriction operators $\mathbf{P}_s$ and $\mathbf{R}_s$ defined in Section \ref{subsec:spatial_discretization}. Note that these operators have to be computed only once and can be used for all MGRIT variants.  Furthermore, nested-iterations, described in Section \ref{subsec:mgrit-alg}, is used to generate an improved initial guess. For this preprocessing step, we use the reference signal \eqref{eq:machine_reference_signal} as a voltage source based on the idea in \cite{Gander_2019aa}. For the iterations of the MGRIT algorithm, the discontinuous signal \eqref{eq:machine_pwm} is used on every grid level. Since we are interested in the algorithmic behavior of different MGRIT variants, in the following, we report on tests with MGRIT convergence measured based on the absolute space-time residual norm instead of the Joule losses \eqref{eq:losses}. More precisely, the convergence tolerance is chosen to be $10^{-4}$, measured in the discrete $L_2$-norm. We have checked that for all algorithms, this fixed space-time residual norm is sufficient to achieve discretization error for our linear model problem if the algorithms converge. That is, algorithms either converge to the time-stepping solution of the linear model of the four-pole squirrel cage induction machine with a PWM voltage source or convergence is only observed after $N_t/m$ or $N_t/(2m)$ iterations for $F$- or $FCF$-relaxation, respectively.

\subsection{Effects of spatial coarsening}
\label{sec:case_study_spatial_coarsening}

We apply several MGRIT variants to the linear machine model problem on the space-time domain $\Omega\times [0,0.03125]$ s. The time interval, corresponding to about one and a half periods, consists of $2^{15}+1$ time steps of size $\Delta t = 2^{-20}$ that are distributed evenly among 256 processors. Note that this time-step size resolves physical processes as well as the pulses of the PWM voltage source. \revision{As we study the initial behavior of the machine, also note that the final time is not sufficient to reach the steady state.} With this choice of number of processors and time points, the temporal coarsening factor is chosen to be $m = 256$ on the first level, resulting in a second-level time grid with $256$ points. The five-level algorithms use a coarsening factor of four on all coarse levels. 

\begin{table*}
	\renewcommand{\arraystretch}{1.3}
	\begin{center}
		\begin{tabular}{ |l|c|c|c|c| }
		\hline
			MGRIT variant  & Iterations & Setup time & Solve time & Total time \\ \Xhline{2\arrayrulewidth}
			two-level cycles with $F$-relax.           & $12$ & $9523$ s & $166{,}265$ s & $175{,}788$ s  \\ 
 			two-level cycles with $FCF$-relax.      & $7$ & $9515$ s & $125{,}149$ s & $134{,}664$ s \\ \hline
 			five-level V-cycles with $FCF$-relax.  & $8$ & $3293$ s & $77{,}146$ s & $80{,}439$ s \\ \hline
			five-level F-cycles with $F$-relax.       & $12$ & $2396$ s & $85{,}117$ s & $87{,}513$ s\\ 
 			five-level F-cycles with $FCF$-relax.  & $7$ & $3291$ s & $87{,}660$ s & $90{,}951$ s\\ \hline
	\end{tabular}
	\end{center}
	\caption{Results for the linear induction machine discretized on a space-time grid of size $69{,}384 \times 2^{15}$  in terms of iterations, setup, and solve time on 256 processors for the different MGRIT variants without spatial coarsening.  
}
	\label{tab:lin_sol_without_sc}
\end{table*}

Table \ref{tab:lin_sol_without_sc} shows the number of iterations, setup, solve, and total time in seconds for the different MGRIT variants without spatial coarsening. The setup time consists of the setup of the algorithm and the generation of an improved initial guess using nested iterations. Note that for nested iterations, the same relaxation scheme as the algorithm itself is used, e.\,g., the setup for the five-level $F$-cycle with $F$-relaxation is faster than that for the five-level $F$-cycle with $FCF$-relaxation. The two-level algorithm with $F$-relaxation reaches the desired tolerance after $12$ MGRIT iterations. While the use of $FCF$-relaxation is beneficial in the two-level setting, reducing both the number of iterations and the runtime, for $F$-cycles, $F$-relaxation performs better than $FCF$-relaxation. Comparing the total runtimes, the multilevel variants are faster than the two-level methods, with the five-level $V$-cycle being the fastest. The problem on the second level is still very expensive, and, therefore, using the two-level method recursively for this problem allows for better performance. 

Now, we consider the same MGRIT variants with spatial coarsening. More precisely, we add spatial coarsening by using the grids $\mesh_2$ and $\mesh_3$, following either the direct or delayed coarsening strategy. Note that in the five-level MGRIT variants, direct spatial coarsening refers to using $\mesh_1$ on the finest grid, $\mesh_2$ on the second level, and $\mesh_3$ on all remaining coarse levels. Delayed spatial coarsening, in contrast, uses $\mesh_3$ only on the coarsest time grid, $\mesh_2$ on the second coarsest level, and $\mesh_1$ on the three finest levels. In the two-level methods, $\mesh_2$ or $\mesh_3$ is used on the coarse grid. Table \ref{tab:lin_sol} demonstrates the effects of adding spatial coarsening to the MGRIT variants considered in Table \ref{tab:lin_sol_without_sc}. All variants with $FCF$-relaxation converge to the solution in the same number of iterations as the variants without spatial coarsening. Especially for this problem, the direct spatial coarsening strategy does not degrade the convergence behavior. The situation is different for all variants with $F$-relaxation. Neither the two-level method with spatial coarsening converges to the desired tolerance in significantly fewer than $N_t/m$ iterations, nor does the five-level $F$-cycle. 

\begin{table*}
	\renewcommand{\arraystretch}{1.3}
	\begin{center}
		\begin{tabular}{ |l|>{\centering}m{4em}|c|c|c|c|>{\centering\arraybackslash}m{5.5em}| }
		\hline
		MGRIT variant & SC\\ strategy & Iterations & Setup time & Solve time & Total time & Speedup\newline w.r.t. no SC\\ \Xhline{2\arrayrulewidth}
		two-level cycles with $F$-relax.                       & direct  &  \xmark & - & - & - & - \\ 
 		two-level cycles with $FCF$-relax. $\left(\Omega_2\right)$  			& direct  &  $7$ & $2444$ s & $73{,}713$ s & $76{,}157$ s & $1.77$ \\
 		 two-level cycles with $FCF$-relax. $\left(\Omega_3\right)$   			& direct  &  $7$ & $685$ s & $48{,}866$ s & $49{,}551$ s & $2.72$ \\ \hline
 		\multirow{2}{*}{five-level V-cycles with $FCF$-relax.} & delayed &  $8$ & $1664$ s & $73{,}022$ s & $74{,}686$ s & $1.08$ \\ 
 					 					 & direct  & $8$ & $316$ s & $66{,}233$ s & $66{,}549$ s & $1.21$ \\ \hline
		\multirow{2}{*}{five-level F-cycles with $F$-relax.} & delayed & \xmark & - &- &- & -\\ 
 					                                                     & direct & \xmark  & - & - &- & -  \\  \hline
 		\multirow{2}{*}{five-level F-cycles with $FCF$-relax. } & delayed &  $7$ & $1665$ s & $74{,}216$ s & $75{,}881$ s & $1.20$ \\ 
 		                                                                                   & direct  & $7$ & $317$ s & $60{,}372$ s & $60{,}689$ s & $1.50$ \\  \hline
	\end{tabular}
	\end{center}
	\caption{Results similar to those of Table \ref{tab:lin_sol_without_sc} but with spatial coarsening; \xmark$\;$indicates no convergence to the desired tolerance in less than $N_t/m=256$ iterations and speedup is measured relative to the same MGRIT variant without spatial coarsening.}
	\label{tab:lin_sol}
\end{table*}

\begin{figure}
	\centering
	\includegraphics[width=0.4\textwidth]{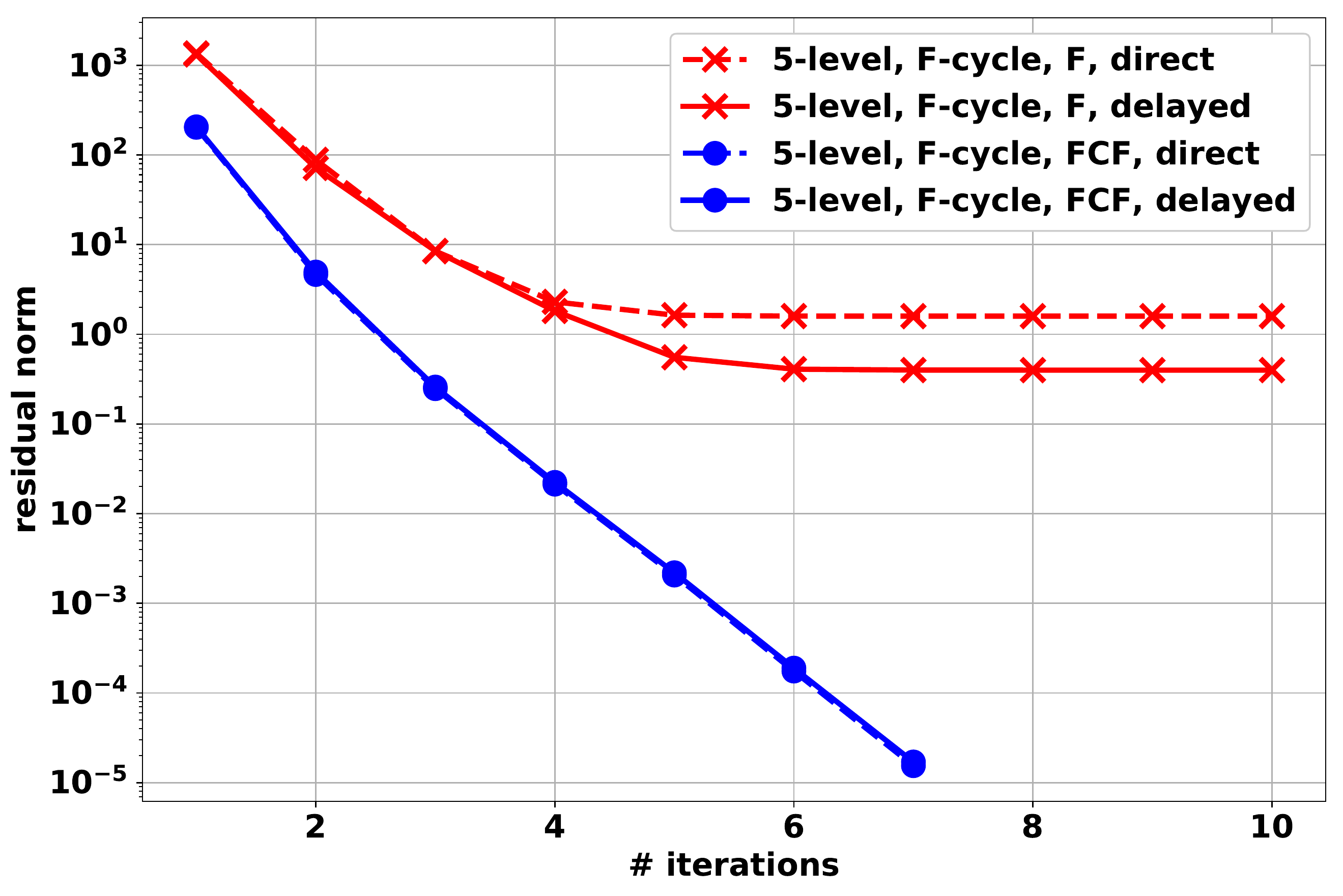}
		\caption{Convergence results for the different five-level $F$-cycle variants with direct (dashed lines) and delayed (solid lines) spatial coarsening.}
	\label{fig:res_norm_with_sc}
\end{figure}

To better understand the degradation in convergence when using $F$-relaxation, Figure \ref{fig:res_norm_with_sc} details the convergence behavior for all $F$-cycle variants considered in Table \ref{tab:lin_sol}. Shown are the space-time residual norms over the first ten MGRIT iterations for all $F$-cycle variants. Dashed lines are results for using direct spatial coarsening and solid lines represent applying the delayed approach. The variants with $FCF$-relaxation show linear convergence behavior, and there is effectively no difference between both spatial coarsening strategies. In contrast, convergence of MGRIT with $F$-relaxation stagnates after five iterations. We assume this behavior comes from the scalar-valued unknowns of our system, which are injected between the spatial grids. The additional $CF$-relaxation in the $FCF$-relaxation scheme improves the spatial interpolation and fixes this problem. 

Comparing the runtimes of the different MGRIT variants with and without spatial coarsening in Tables \ref{tab:lin_sol_without_sc} and \ref{tab:lin_sol}, the two-level algorithms show the best speedup, decreasing the total runtime by up to a factor of about $2.72$, followed by the five-level F-cycle with $FCF$-relaxation and direct spatial coarsening, which is faster than the same variant without spatial coarsening by a factor of about $1.5$. Note that, for this problem, the use of direct spatial coarsening does not degrade MGRIT convergence and, thus, choosing this more aggressive coarsening strategy leads to higher speedups than applying delayed spatial coarsening. Applying a more aggressive spatial coarsening strategy also in the two-level method with $FCF$-relaxation similarly allows for a larger speedup over the same algorithm without spatial coarsening. Using the spatial mesh $\mesh_3$ instead of $\mesh_2$ on the coarse time grid, reduces the runtime from $76{,}157$ seconds to $49{,}551$ seconds and, thus, increases the speedup from a factor of about 1.77 to a factor of about 2.72. However, the scalability of the two-level method is limited by the size of the coarse-grid problem, and spatial coarsening can lead to divergence in the nonlinear setting, as will be seen in Section \ref{subsec:nonlinear_effecs_spatial_coarsening}.

\subsection{Discussion}
\label{sec:lin_discussion}

Summarizing, the MGRIT algorithm is applied successfully to the linear material model of an induction machine with a pulsed voltage source. With our choice of spatial grid transfer operators, spatial coarsening can be added for all variants with $FCF$-relaxation, while the methods with $F$-relaxation methods do not convergence to the time-stepping solution in less than $N_t/m$ iterations. A better choice could resolve this behavior. For our problem, we see no difference in the direct and delayed strategy, and, therefore, the direct approach should be used for reducing the runtime. The speedup by adding spatial coarsening corresponds directly to the amount of work on the coarser grids, i.\,e., the two-level algorithm has the best speedup, followed by the $F$-cycle.

\section{Full nonlinear model}
\label{sec:nonlinear_model}

In this section, we use the results from Section \ref{sec:case_study_spatial_coarsening} of different MGRIT variants compared with one another for studying the application of MGRIT to the full nonlinear model of the induction machine ``im\_$3$\_kw'', i.\,e., with a state dependent reluctivity $\nu$ which models magnetic saturation. In particular, we are interested in effects of spatial coarsening in this setting as well as in the parallel performance of MGRIT. Furthermore, we compare MGRIT with sequential time stepping.

\subsection{Numerical parameters}
\label{nonlinear_numerical_parameters}

To limit computational costs in the full nonlinear setting, we perform experiments on a smaller space-time domain. More precisly, we consider $\mesh_2$ as the spatial grid and choose $\Delta t = 2^{-20}$ and $N_t=10{,}753$ time steps, resulting in a final time $t_f \approx 0.01$ s. Convergence of MGRIT is measured based on the relative change in values of joule losses at C-points of two consecutive iterations. MGRIT iterations are stopped when the maximum norm of the relative difference of two successive iterates at C-points is below $10^{-2}$, i.\,e., when the maximum relative change in joule losses is smaller than 1\% for all C-points. Note that this stopping criterion is much looser compared to that in Section \ref{sec:linear_model}, but it corresponds to a choice made for prototyping of induction machines in a realistic setting. However, in contrast to using a convergence tolerance based on the space-time residual, it is not guaranteed that MGRIT converges to the discrete time-stepping solution. For the MGRIT variants considered in this section, we have checked that, if convergence is observed, iterates indeed converge to the discrete time-stepping solution.

In the full nonlinear setting, each application of the time propagator $\boldsymbol{\Phi}$ involves an iterative nonlinear solve. GetDP, called by the time stepper routine of our\linebreak MGRIT implementation, uses Newton's method with damping for the solution of the nonlinear problems. As convergence of the Newton solver can depend on the initial guess, in Section \ref{subsec:nonlinear_effecs_spatial_coarsening}, we consider not only the effect of spatial coarsening on MGRIT convergence, but also on the convergence of the Newton solver. We then use these results in Section \ref{sub:parallel_results} to choose a set of MGRIT variants for a strong scaling study and comparison to sequential time stepping.

\subsection{Effects of spatial coarsening}
\label{subsec:nonlinear_effecs_spatial_coarsening}

The results in Section \ref{sec:case_study_spatial_coarsening} show that adding spatial coarsening in MGRIT is beneficial in the case of $FCF$-relaxation, whereas MGRIT with $F$-relaxation does not converge when spatial coarsening is added. We now consider two- and five-level MGRIT variants with $FCF$-relaxation and compare the performance of these methods without and with direct spatial coarsening using the spatial mesh $\Omega_3$ on coarse temporal grids. Additionally, we look at the performance of two-level MGRIT with $F$-relaxation and without spatial coarsening. All tests are performed using 256 processors. Choosing again the non-uniform temporal coarsening strategy described in Section \ref{subsec:linear_algorithmic_parameters}, a coarsening factor of $42$ is used on the first level, and factor-four coarsening is applied on all coarse levels. 

Table \ref{tab:results_nonlin} shows iteration counts and runtimes of the different MGRIT variants, with timings split up again into setup and solve times. No results are given for the two-level method with $FCF$-relaxation and with spatial coarsening, since during the $F$-relaxation step of the first MGRIT iteration, the Newton solver does not converge for at least one time step and, thus, the algorithm cannot be applied in this setting. In contrast to the setting of linear material laws, here, an initial guess that is in the region of convergence of the Newton solver is crucial. Applying a nested iterations strategy in the two-level setting, the initial guess at the $C$-points of the fine  grid are given by the interpolated solutions of the coarse grid. As no values are given at $F$-points, $F$-relaxation is then carried out using the solution from the previous time point as the initial guess, i.\,e., for the first $F$-point in each $F$-interval, the solution at the preceding $C$-point is used as the initial guess. After solving the corresponding nonlinear problem at the first $F$-point in each interval, the solution is taken as the initial guess for the following $F$-point, and so on. However, this strategy appears to result in a poor initial guess for the Newton solver at some time points. There are several difficulties for obtaining a good initial guess, including the use of spatial coarsening, the discontinuities of the right-hand side of the problem, as well as the different time-step sizes on the temporal grids in the multigrid hierarchy. Improving the initial guess requires further investigation of strategies for tackling these challenges. This is a topic of future work.

\begin{table*}
	\renewcommand{\arraystretch}{1.3}
	\begin{center}
		\begin{tabular}{ |l|>{\centering}m{4em}|c|c|c|c| }
		\hline
		MGRIT variant & Spatial\\ coarsening & Iterations & Setup time & Solve time & Total time\\ \Xhline{2\arrayrulewidth}
		two-level cycles with $F$-relax.  & no  & $4$ & $15{,}054$ s & $74{,}150$ s & $89{,}204$ s \\ 
 		two-level cycles with $FCF$-relax. & yes  & \xmark & - & - & - \\ \hline
 		\multirow{2}{*}{five-level V-cycles with $FCF$-relax.} & no & $3$ & $5870$ s & $25{,}509$ s & $31{,}379$ s\\ 
 		                                                                                       & yes & $3$ & $1297$ s & $18{,}763$ s & $20{,}060$ s\\ \hline
		\multirow{2}{*}{five-level F-cycles with $FCF$-relax.} & no & $3$ & $5861$ s& $41{,}138$ s & $46{,}999$ s \\ 
		                                                                                  & yes & $3$ &  $1301$ s & $22{,}430$ s & $23{,}731$ s \\ \hline
		\end{tabular}
		
	\end{center}
	\caption{Number of iterations, setup, solve and total time on 256 processors of various MGRIT variants applied to the full nonlinear model of the induction machine ``im\_$3$\_kw'' discretized on a space-time grid of size $17{,}496 \times 10{,}753$. \xmark$\;$indicates no convergence of at least one nonlinear spatial solve within GetDP. }
	\label{tab:results_nonlin}
\end{table*}

Looking at the total runtimes of the different variants without spatial coarsening, the five-level $V$-cycle algorithm is fastest, followed by the five-level $F$-cycle variant which is about a factor of 1.5 times slower than considering $V$-cycles. Furthermore, five-level $F$-cycles are already about twice as fast as the two-level method. Adding spatial coarsening in the multilevel schemes, we can benefit over the two-level algorithm even more. Considering five-level $V$- and $F$-cycles with spatial coarsening, the factor in comparison with the runtime of the two-level method can be increased from 1.9 or 2.8 when applying $F$- or $V$-cycles without spatial coarsening, respectively, to a factor of about 3.8 or 4.4 when spatial coarsening is added.

\subsection{Parallel results}
\label{sub:parallel_results}

Now we present strong scaling results of the five convergent MGRIT variants considered in Table \ref{tab:results_nonlin}. We perform parallel tests using between eight and $256$ processors for parallelization in time. The temporal coarsening strategy in all runs is fixed at applying a coarsening factor of $42$ on the finest grid and factor-four coarsening on all coarse grids. Additionally, we are interested in the benefits over sequential time stepping.  

\begin{figure*}
	\begin{center}
	\includegraphics[width=0.7\textwidth]{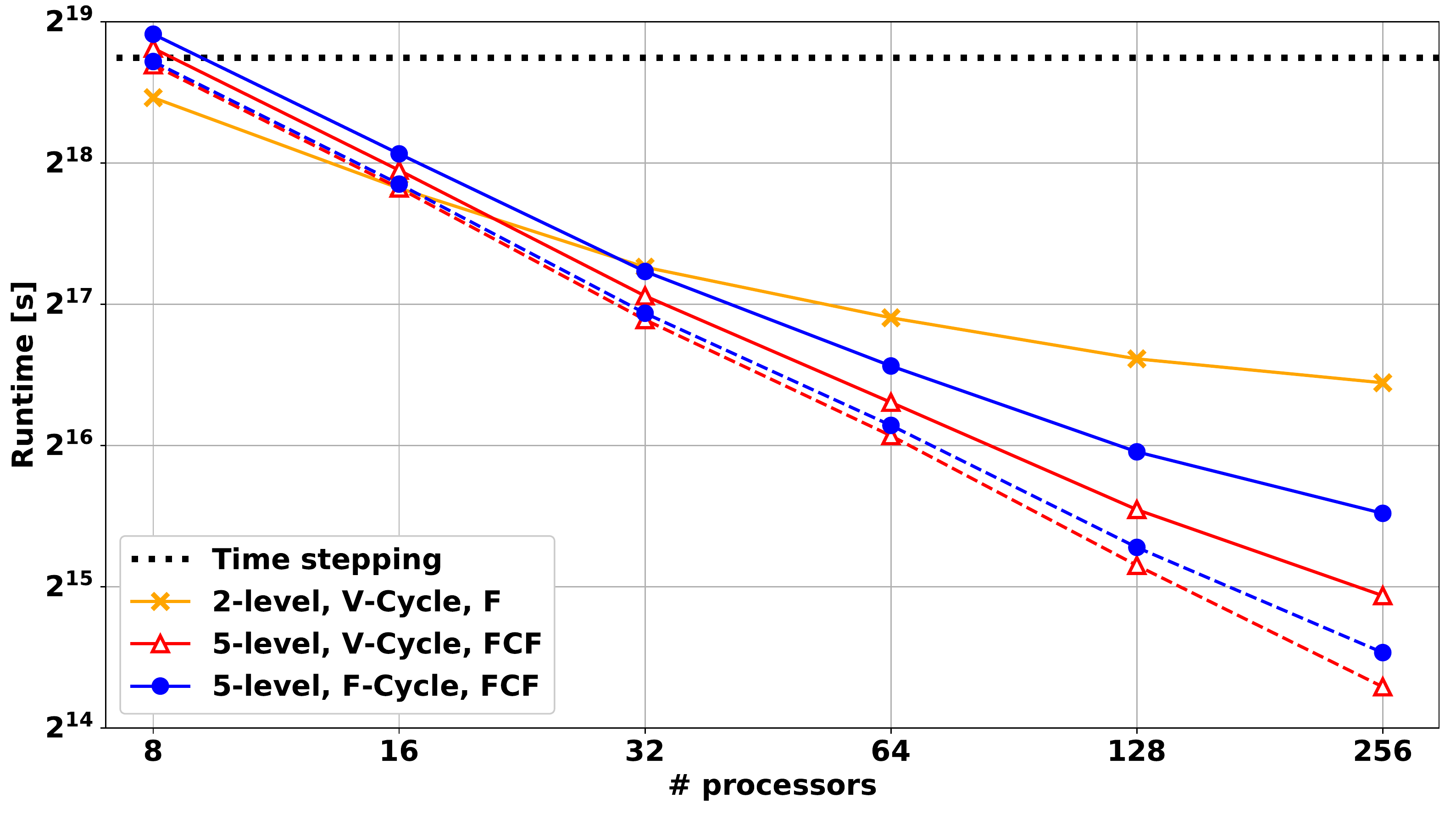}
		\caption{Total time-to-solution for the nonlinear induction machine ``im\_$3$\_kw'' using different MGRIT variants and sequential time stepping. Solid lines are runtimes without spatial coarsening, and dashed lines represent results with spatial coarsening. The dotted line shows the runtime of time stepping on one processor for reference purposes.}
	\label{fig:strong_scaling}
	\end{center}
\end{figure*}

Figure \ref{fig:strong_scaling} shows total runtimes of the different MGRIT variants as a function of the number of processors, as well as the time-to-solution of the sequential block forward solve for reference purposes. For the multilevel variants, runtimes are shown for using spatial coarsening (dashed lines) and without spatial coarsening (solid lines). The runtime of all multilevel variants using eight processors as well as the time-to-solution of sequential time stepping is about five to six days, while using eight-way parallelism in the two-level method results in a runtime of only about four days. However, the multilevel variants show better strong parallel scaling and, thus, lead to faster time-to-solution at large processor counts than the two-level method. Increasing the number of processors to $256$, the total runtime can be reduced to about $8.72$ hours or $5.57$ hours, when considering MGRIT $V$-cycles without or with spatial coarsening, respectively. While $F$-cycles with spatial coarsening allow for a similar reduction, $F$-cycles without spatial coarsening using $256$-way parallelism already take about half a the two-level variant reduces the runtime to about one day.

Table \ref{tab:efficiency_and_speedup} details speedups and parallel efficiencies for the MGRIT variants. The speedup is computed relative to sequential time stepping, representing the algorithm with the minimum runtime on one processor, i.\,e., the speedup using $p$ processors is given by $S_p=T_\text{TS}/T_\text{MGRIT}(p)$. The parallel efficiency is measured as $S_p/p$ for $p$ processors. Using multilevel methods, we can benefit more over sequential time stepping at larger processor counts than with the two-level algorithm. For example, while on $16$ processors the speedup is similar for all methods, using $32$ processors results in a speedup of up to a factor of $3.6$ for a five-level variant, whereas the two-level method yields only a speedup of a factor of $2.8$. Increasing the number of processors to $256$, with the two-level method we achieve a speedup of a factor of about $5$. In contrast, multilevel variants at least double the speedup factor. Considering spatial coarsening, leads to the best speedup with the $V$-cycle algorithm being nearly $22$ times faster than the sequential time stepping method. Moreover, the use of spatial coarsening in MGRIT $F$-cycle is crucial for improving the parallel scalability for more than $64$ processors, since this allows reducing computations on coarse levels that dominate over communication costs. While speedups and efficiencies degrade for $F$-cycles without spatial coarsening, the degredation is only modest when adding spatial coarsening. Note that we bind the temporal coarsening strategy to the maximum number of available processors and, thus, faster runtimes on eight to 128 processors may be possible by adjusting the temporal coarsening strategy to the number of processors.

\begin{table*}
	\begin{center}
		\begin{tabular}{|C{5em}|C{4em}|C{4em}|C{4em}|C{4em}|C{4em}|C{4em}|}
		\hline 
		$p$ & 8 & 16 & 32 & 64 & 126 & 256 \\ \hline
		\multicolumn{7}{ |c| }{two-level cycles with $F$-relax} \\ \hline
		Speedup & 1.22 & 1.90 & 2.80 & 3.59 & 4.39 & 4.93 \\
		Efficiency & 15.22\% & 11.85\% & 8.74\% & 5.60\% & 3.43\% & 1.93\% \\ \hline 
		\multicolumn{7}{ |c| }{five-level $V$-cycle with $FCF$-relax.} \\ \hline
		Speedup & 0.96 & 1.74 & 3.22 & 5.43 & 9.19 & 14.02 \\
		Efficiency & 11.94\% & 10.86\% & 10.07\% & 8.48\% & 7.18\% & 5.48\% \\ \hline 
		\multicolumn{7}{ |c| }{five-level $F$-cycle with $FCF$-relax.} \\ \hline
		Speedup & 0.89 & 1.60 & 2.86 & 4.54 & 6.92 & 9.36 \\
		Efficiency & 11.14\% & 10.02\% & 8.93\% & 7.10\% & 5.41\% & 3.66\% \\ \hline 
		\multicolumn{7}{ |c| }{five-level $V$-cycle with $FCF$-relax. and SC} \\ \hline
		Speedup & 1.04 & 1.90 & 3.62 & 6.40 & 12.01 & 21.93 \\
		Efficiency & 12.96\% & 11.86\% & 11.31\% & 9.99\% & 9.44\% & 8.57\% \\ \hline 
		\multicolumn{7}{ |c| }{five-level $F$-cycle with $FCF$-relax. and SC} \\ \hline
		Speedup & 1.02 & 1.86 & 3.51 & 6.08 & 11.06 & 18.54 \\
		Efficiency & 12.74\% & 11.63\% & 10.96\% & 9.50\% & 8.64\% & 7.24\% \\ \hline 
		
\end{tabular}
\label{tab:efficiency_and_speedup}
\end{center}
\caption{Speedup and efficiency of different MGRIT variants using various number of processors. The speedup is given relative to the time-to-solution for the sequential time-stepping on one processor. Parallel efficiency is measured as $S_p/p$, where $S_p$ is the speedup for $p$ processors.}\label{tab:efficiency_and_speedup}
\end{table*}

\section{Conclusions}
\label{sec:conclusions}
In this paper, the MGRIT algorithm is applied to a model of a realistic two-dimensional electrical machine with a pulsed excitation. Using the non-intrusive character of MGRIT for an existing model of an induction machine, performance of MGRIT is dominated by the cost of the time-stepping routine which carries out a nonlinear spatial solve for a given time step. To reduce the cost of spatial solves on coarse grids, the use of spatial coarsening is investigated. We demonstrate that spatial coarsening is a powerful option that can significantly reduce the runtime of the algorithm, but it also comes along with new challenges and problems, e.\,g., degradation in MGRIT convergence, as also observed in \cite{MR3716560} for the $p$-Laplacian, and divergence of the Newton solver due to a poor initial guess. These are topics of future work.

Strong scaling results show significant speedups compared to sequential time stepping. For a full nonlinear model of an induction machine and using moderate numbers of processors, we demonstrate that MGRIT allows reducing the simulation time from several days to only a few hours. While MGRIT is a non-intrusive approach, there is still a very large parameter space for the algorithm, such as the choice of the level structure, temporal and spatial coarsening strategies, and so on. Especially if computational resources are limited to small or moderate numbers of processors, good parameter choices are crucial for obtaining the best possible speedup over sequential time stepping. Future work will focus on developing a model for optimizing MGRIT performance for given numbers of processors and time points.

\begin{acknowledgements}
 The work is supported by the Excellence Initiative of the German Federal and State Governments, the Graduate School of Computational Engineering at TU Darmstadt, and the BMBF (project PASIROM; grants 05M18RDA and 05M18PXB). The authors would like to thank Iryna Kulchytska-Ruchka for providing code for GetDP simulations of the ``im\_3\_kW'' model. 
\end{acknowledgements}

%
%

\bibliographystyle{spmpsci}      
\bibliography{pint_refs}   

%
%

\end{document}